\documentclass[12pt]{article}

\def\({\left(}
\def\){\right)}

\def\bk{\bigskip }

\def\sk{\smallskip }

\def\n{\noindent }

\def\barr{\begin{array}}
\def\earr{\end{array}}

\def\bit{\begin{itemize}}

\def\eit{\end{itemize}}

\usepackage{amsmath, amsthm, amscd, amsfonts, amssymb}
\usepackage[usenames]{color}

\usepackage{mathrsfs}

\numberwithin{equation}{section} 

\newtheorem{theorem}{Theorem}[section]
\newtheorem{proposition}[theorem]{Proposition}

\newtheorem{lemma}[theorem]{Lemma}
\theoremstyle{definition}

\newtheorem{remark}[theorem]{Remark}

\def\calo{\mathcal{O}}

\def\bbe{{\mathbb{E}}}

\def\bbr{{\mathbb{R}}}
\def\bbp{{\mathbb{P}}}

\def\1{^{-1}}

\def\calo{{\mathcal{O}}}

\def\9{{\infty}}

\def\a{{\alpha}}

\def\g{{\gamma}}
\def\wt{\widetilde}

\def\vf{{\varphi}}

\def\barr{\begin{array}}
\def\earr{\end{array}}

\def\bk{\bigskip }
\def\sk{\smallskip}
\def\n{\noindent }

\def\({\left(}
\def\){\right)}
\def\<{\left<}
\def\>{\right>}

\def\wt{\widetilde}
\def\wh{\widehat}

\def\ve{{\varepsilon}}

\def\e{{\epsilon}}

\begin{document}

\begin{center}
{\Large{\bf Gaussian fluctuations and moderate deviations of eigenvalues in
unitary invariant ensembles}}
\bigskip\bk

{\large{\bf Deng
Zhang}}\footnote{School of Mathematical Sciences, Shanghai Jiao Tong
University, 200240 Shanghai, China. E-mail address: dzhang@sjtu.edu.cn}
\end{center}

\bk\bk\bk

\begin{quote}
\n{\small{\bf Abstract.} We study the limiting behavior of the $k$-th eigenvalue $x_k$ of
unitary invariant ensembles
with Freud-type and  uniform convex potentials.
As both $k$
and $n-k$ tend to infinity, we obtain Gaussian fluctuations for
$x_k$ in the bulk and soft edge cases, respectively.
Multi-dimensional central limit
theorems, as well as moderate deviations, are also proved.
This work generalizes earlier results in the GUE  and unitary invariant
ensembles with monomial potentials of even degree.
In particular, we obtain the precise asymptotics of
corresponding Christoffel-Darboux kernels as well. } \\

\end{quote}

\section{Introduction and main results} \label{Intro}
We are concerned with the unitary invariant ensemble of $n \times n$ Hermitian
matrices $\mathscr{H}_n$ with the probability distribution defined by
\begin{align} \label{UE}
    \mathbb{P}_n(dH)=C_ne^{-nTrV(H)}dH,\ \ H\in \mathscr{H}_n,
\end{align}
where $C_n$ is a normalization constant, $V(x)$ is an external
potential, which is real analytic and satisfies
$V(x)/\log(x^2+1)\rightarrow\infty$, as $|x|\rightarrow\infty$, and
$dH$ stands for the Lebesgue measure on the algebraically
independent entries of $H$, i.e., $
    dH=\prod_{1\leq i<j\leq n} d{\rm Re} H_{ij}d {\rm Im}H_{ij} \prod_{i=1}^n
    dH_{ii}.$

It is well known (cf.\cite{D99}) that the
distribution \eqref{UE} induces a probability density function of
the corresponding $n$ ordered real eigenvalues $\{x_i\}_{i=1}^n$,
$x_1<...<x_n$, given by
\begin{eqnarray} \label{UE-density}
    \mathscr{R}_{n,n}(x_1,...,x_n)=\frac{1}{Z_n}\prod_{1\leq i<j\leq
    n}|x_i-x_j|^2\ \exp\(-n\ \sum\limits_{i=1}^nV(x_i)\),
\end{eqnarray}
where $Z_n$ is a normalization. In particular,
the quadratic potential (i.e.,
$V(x)=2x^2$) corresponds to the classical Gaussian Unitary Ensemble (GUE).

Unitary invariant ensembles have been extensively studied in literature.
The physical significance of  probability
distribution \eqref{UE-density} can be interpreted as a Gibbs
measure for $n$ identical charged particles in $\mathbb{R}$, at the
inverse temperature $\beta= 2$,  with a logarithmic interaction and
with an external potential $V$. The generalization to general
inverse temperatures $\beta>0$ is known as beta ensembles or
log-gases. See, e.g., \cite{E12,F10}.

One remarkable global
property is that the $1$-point correlation function of
\eqref{UE-density} converges weakly to an equilibrium measure,
which is the well-known
semicircle law in the GUE (see \cite{D99,J98}).
Moreover, the dynamical interpretations of \eqref{UE-density} and
equilibrium measure are related closely to the generalized
Dyson Brownian motion and Mackean-Vlasov equation respectively. We
refer, e.g., to  \cite{C92, RS93} for the GUE case and the
recent work \cite{LLX13} for  beta ensembles.

The main interests of this paper are concerned with  local
fluctuations, as well as  moderate deviations,
of the $k$-th eigenvalue $x_k$ of
a matrix taken randomly from a unitary ensemble,
when both
$k$ and $n-k$ tend to infinity.
Local fluctuations of the $k$-th eigenvalue $x_k$, in a general context,
turn out to be universal.

On the one hand, when $k$ or $n-k$ is fixed,
these fluctuations obey
the celebrated Tracy-Widom distribution.
We refer to \cite{TW94} for the GUE,
\cite{DG07} and \cite{PS03} for unitary invariant ensembles with
Freud-type potentials and
analytic potentials, respectively.
See also \cite{BEY14} for general beta ensembles and \cite{S09}
for orthogonal ensembles.

On the other hand, when both $k$ and $n-k$ tend to infinity, the $k$-th eigenvalue $x_k$
is asymptotically normally distributed in both the bulk and
edge cases. This result was first
proved by Gustavsson \cite{G05} for the GUE and later extended to various other matrix models.
In the bulk case,
see \cite{Z14} for unitary invariant ensembles with monomial potentials,
and  \cite{L16} for more general potentials  including  uniform convex potentials.
In the  edge case,
we refer to \cite{BEY14} for  general beta
ensembles with   potentials independent of $n$.
See also \cite{S06} for complex covariance matrices,  \cite{TV11} for Wigner Hermitian
matrices, and \cite{O09}
for real symmetric Wigner matrices.

Moreover, it is well known (see e.g. \cite{DGZ03, DE13.0, Z16})
that moderate deviations are related closely to central limit theorems.
Recently,
a moderate deviation principle was obtained
for general determinantal point processes in \cite[Theorem 1.4]{DE13}.
This result implies the moderate deviations of  $k$-th eigenvalue of
Wigner matrix
in the bulk and edge cases (\cite{DE13, DE13.2}),
which is indeed another motivation of the present work.

Here we are mainly concerned with
unitary ensembles with Freud-type weights
or uniform convex potentials, that is,

$(i)$ Freud-type potential,
\begin{equation} \label{v.2}
   (V(x):=) V_n(x)=\frac{1}{n}Q(c_nx+d_n),
\end{equation}
where $Q(x)=\sum_{k=0}^{2m}q_kx^k,m\in\mathbb{N}^+$,
$q_{2m}=\frac{\Gamma(m)\Gamma(\frac{1}{2})}{\Gamma(\frac{2m+1}{2})}$,
$m\geq 1$,
$c_n=\frac 12 (\beta_n-\alpha_n)$, $d_n=\frac 12 (\beta_n+\alpha_n)$, and
$\alpha_n$ $\beta_n$ are the n-th Mhasker-Rakhmanov-Saff numbers
(see Section \ref{RHM} below for details).

$(ii)$ Uniform convex potential,
\begin{equation} \label{v.3}
    \inf\limits_{x\in \mathbb{R}}V''(x) \geq c >0
\end{equation}
for some $c>0$.

We obtain the Gaussian fluctuations of $k$-th eigenvalue  in both the bulk and edge cases.
Multi-dimensional central limit theorems of eigenvalues are obtained as well.
These generalize earlier results in the GUE
\cite{G05} and  unitary invariant ensembles with
monomial potentials \cite{Z14}.
Furthermore,
the moderate deviations of $k$-th eigenvalue are obtained in both the bulk and edge cases as well,
thereby generalizing the  results in the GUE  \cite{DE13, DE13.2}
to unitary invariant ensembles with Freud-type and uniform convex potentials.

The proof is mainly
based on the central limit theorems in \cite{CL95,S00,S00.2}
and the moderate deviation principle in \cite{DE13}.
In particular,
we apply the Riemann-Hilbert approach,
developed in \cite{DKMVZ99.0,DKMVZ99},
to obtain the precise asymptotics of corresponding Christoffel-Darboux kernels,
which enable us to obtain the asymptotics of
expectation and variance for the counting statistic for an interval
and would be also of independent interest.

After this paper was finished,
we learned about the works \cite{BEY14} and \cite{L16}
where the Gaussian fluctuations
are proved in the edge and bulk cases, respectively,
for uniform convex potentials
(and also other general potentials or matrix models)
based on different approaches.
We would like to refer to \cite{BEY14, L16} for more details.

Below we formulate the main results of this paper.
Recall that the equilibrium measure $\mu_{V}$ is the unique
minimizer of  variational problem
\begin{equation} \label{vp.1}
    \mu_{V}={\arg min}_{\mu\in\mathscr{M}_1(\mathbb{R})} I_{V}(\mu),
\end{equation}
where $\mathscr{M}_1(\mathbb{R})=\{\mu:\int_{\mathbb{R}}d\mu=1\}$,
$I_{V}$ is the Voiculescu free entropy defined by
\begin{equation} \label{vp.2}
    I_{V}(\mu)=\iint \log|s-t|^{-1}d\mu(s)d\mu(t)+\int V(t)d\mu(t).
\end{equation}
For Freud-type and uniform convex potentials,
it is known (see \cite[(4.17)]{DKMVZ99}, \cite[(1.4),(1.5)]{DKMVZ99.0}) that
$\mu_V$ has the density function $\rho_V$ supported on $[b,a]$.
With suitable scaling, we may assume that $b=-1$
and $a=1$ without lose of generality.
The density function $\rho_{V}$  can be  also characterized by the Euler-Lagrange
equations below (cf. \cite[(4.18), (4.19)]{DKMVZ99}, \cite[(1.10), (1.11)]{DKMVZ99.0})
\begin{align}
  2 \int \log |x-s|  \rho_{V}(s) ds -V(x) &= l,\ x\in[-1,1], \label{rho-vl1} \\
  2 \int \log |x-s| \rho_{V}(s) ds -V(x) &\leq l,
  ~x\in\mathbb{R}/[-1,1].  \label{rho-vl2}
\end{align}

The Gaussian fluctuation results
are formulated in Theorem
\ref{Bulk-Thm} and \ref{Edge-Thm} below.
As in \cite{G05},
we use the notation  $k(n) \sim n^\theta$
to mean that $k(n) = h(n) n^\theta$,
where $h$ is any function satisfying
$h(n) n^{-\epsilon} \to 0$ and
$h(n) n^\epsilon \to \9$
as $n\to \9$ for all $\epsilon >0$.

\begin{theorem}\label{Bulk-Thm} ({\it Bulk case.})
Consider the unitary invariant ensemble \eqref{UE} with the Freud-type and
uniform convex potential as in \eqref{v.2} and \eqref{v.3}, respectively.

$(i)$. Let $G(s) =\int_{-1}^{s}\rho_{V}(x)dx$, $-1\leq s\leq 1$,
and $t=t(k,n)=G^{-1}(k/n)$, where $k=k(n)\in[cn,(1-c)n]$, $c\in(0,
1/2)$. Set
\begin{equation} \label{def-Xn}
   X_n:=  \frac{x_k-t}{\frac{\sqrt{\log n}}{\sqrt{2\pi^2}n\rho_{V}(t)}}.
\end{equation}
Then, $X_n \rightarrow  N(0,1)$ in distribution, as $n\rightarrow \infty$.

$(ii)$. Let $\{x_{k_i}\}_{i=1}^m$ be eigenvalues such that
$0<k_i-k_{i+1}\sim n^{\theta_i},0<\theta_i\leq 1$, and
$k_i\in[c_in,(1-c_i)n]$, $c_i\in(0, 1/2)$. Set
$s_i=s_i(k_i,n)=G^{-1}(k_i/n)$ and
\begin{equation*}
    X_{i,n}:=\frac{x_{k_i}-s_i}{\frac{\sqrt{\log n}}{\sqrt{2\pi^2}n\rho_{V}(s_i)}},\ \ 1\leq i\leq m.
\end{equation*}
Then, for any $\xi_i \in \mathbb{R}$, $1\leq i \leq m$, as $n\rightarrow\infty$,
\begin{equation*}
  \mathbb{P}_n[X_{1,n}\leq \xi_1,...,X_{m,n}\leq
\xi_m]\rightarrow\Phi_{\Lambda}(\xi_1,...,\xi_m).
\end{equation*}
Here $\Phi_{\Lambda}$ is the m-dimensional Normal distribution
function with mean zero and the correlation matrix $\Lambda$,
$\Lambda_{i,i}=1$, $1\leq i\leq m$, and  $\Lambda
_{i,j}=1-\max_{i\leq k<j}\theta_k$, $1\leq i<j\leq m$.
\end{theorem}

\begin{theorem} \label{Edge-Thm} ({\it Edge case.})
Consider the unitary invariant ensemble \eqref{UE} with the Freud-type and
uniform convex potential as in \eqref{v.2} and \eqref{v.3}, respectively.

$(i)$. Let $k$ be such that $k\to \9$ and $k/n \to 0$, as
$n\to \9$. Set
 \begin{align} \label{def-Yn}
    Y_n : =\frac{3\sqrt{2}\pi  a_1 ^{\frac 23}}{2}  \frac{x_{n-k}-[1-(\frac{k}{a_1n})^{\frac{2}{3}}]}
    {\frac{\sqrt{\log k}}{n^{\frac{2}{3}}k^{\frac{1}{3}}}},
 \end{align}
where  for the Freud-type potential  $a_1=\frac{2\sqrt{2}}{3\pi}
  \sum_{i=0}^{m-1}\frac{A_{m-1-i}}{A_m}$ with
$A_j=\prod_{i=1}^j\frac{2i-1}{2i}$,
$A_0=1$, $1\leq j\leq m$,
while for the uniform convex potential
$a_1= \frac{\sqrt{2}}{3\pi} h(1)$ with
$h$ as in Lemma \ref{rhov} below.

Then, $Y_n \to N(0,1)$ in distribution.

$(ii)$. Let $\{x_{k_i}\}_{i=1}^m$ be eigenvalues such that $k_1\sim
n^\gamma$, $0<\gamma<1$, and $0<k_{i+1}-k_i \sim n^{\theta_i}$,
$0<\theta_i<\gamma$.
For the Freud-type potential,
we assume additionally that
$\theta >  \g-\frac{1}{2m}$.
Set
\begin{align*}
    Y_{i,n}:=\frac{3\sqrt{2}\pi a_1^{\frac 23}}{2}  \frac{x_{n-k_i}-[1-(\frac{k_i}{a_1n})^{\frac{2}{3}}]}
    {\frac{\sqrt{\log k_i}}{n^{\frac{2}{3}}k_i^{\frac{1}{3}}}},\ \
    1\leq i\leq m.
 \end{align*}
Then, for any $\xi_i \in \mathbb{R}$, $1\leq i\leq m$, as $n\rightarrow\infty$,
\begin{equation*}
  \mathbb{P}_n[Y_{1,n}\leq \xi_1,...,Y_{m,n}\leq
\xi_m]\rightarrow\Phi_{\Lambda}(\xi_1,...,\xi_m),
\end{equation*}
where $\Lambda$ is as in Theorem \ref{Bulk-Thm}, but with $\Lambda
_{i,j}=1-\g^{-1}\max_{i\leq k<j}\theta_k$,
$1\leq i<j\leq m$.
\end{theorem}

\begin{remark}
For the Freud-type potential,
the additional condition $\theta >  \g-\frac{1}{2m}$
arises in the delicate  estimate of remaining term of $\wt{a}(t)$ in Proposition \ref{Edge-Exp}.
For more details see the proof of \eqref{Dsij-edge} below.
In particular,
when $\g \leq  \frac{1}{2m} $,
one can take any $0<\theta<\g$.
\end{remark}

\begin{remark}
Theorems \ref{Bulk-Thm} and \ref{Edge-Thm} generalize the results in the GUE \cite{G05}
and the unitary invariant ensembles with monomial potentials \cite{Z14}.
After this work was finished,
we learned about the works \cite{BEY14} and \cite{L16}
which also obtained the Gaussian fluctuations
in the edge and bulk cases, respectively,
for the uniform convex potentials and
also for other general potentials or matrix models.
The proof presented below is different, based on the
Riemann-Hilbert approach,
and also gives the precise asymptotics of corresponding
Christoffel-Darboux kernels
(see Section \ref{K-K} and Lemma \ref{5.K.0} below).
\end{remark}

Regarding the moderate deviations, we recall that a sequence
of probability measures $\{\mu_n\} \subseteq \mathscr{M}_1(\bbr)$ is said to satisfy the
large deviation principle with speed $s_n\to \9$ and good rate function
$I: \bbr \to [0,\9]$, if the level sets $\{x\in\bbr: I(x)\leq c\}$ are compact for all
$c\in [0,\9)$ and if for all Borel set $A$ of $\bbr$,
\begin{align*}
  -\inf\limits_{x\in A^o} I(x)
  \leq& \liminf\limits_{n\to\9} \frac{1}{s_n} \log \mu_n (
  A)
  \leq \limsup\limits_{n\to\9} \frac{1}{s_n} \log \mu_n (
  A)
  \leq - \inf\limits_{x\in \overline{A}} I(x),
\end{align*}
where $A^o$ and $ \overline{A}$ denote the interior and closure of
$A$, respectively. In that case, we simply say that $\{\mu_n\}$
satisfies the $LDP(s_n, I)$. We also say that a
family of real valued random variables satisfies the
$LDP(s_n, I)$ if the family of their laws does. In particular, if the deviation scale of random variables
is between that of the law of large number and that of the central limit theorem, this sequence of random variables
is said to satisfy the moderate deviation principle.

We set
$\wt{a}(t) := (1-t)^{-\frac 32} \int_t^1 \rho_{V}(x)dx$
for $t\in [1-\delta, 1]$
with $\delta>0$ small enough.
In particular,
$\wt{a}(t) = a_1(1+o(1))$,
where $a_1$ is the constant as in Theorem \ref{Edge-Thm}.

\begin{theorem} \label{MDP-Thm}
Consider the unitary invariant ensemble \eqref{UE} with the Freud-type and
uniform convex potential as in \eqref{v.2} and \eqref{v.3}, respectively.

$(i)$ ({\it Bulk case.})
Let $k=k(n)\in[cn,(1-c)n]$ with $c\in(0,
1/2)$
and $t=t(k,n)=G^{-1}(k/n)$,
where  $G$ is as in Theorem \ref{Bulk-Thm}.
Let $X_n$ be as in \eqref{def-Xn}.

Then,
for any sequence $\{\g_n\}$ such that $1\ll \g_n \ll \sqrt{\log n}$,
$\{\g^{-1}_nX_n\}$ satisfies the $LDP(\g_n^2, x^2/2)$.

$(ii)$ ({\it Edge case.})
Let $k$ be such that $k\to \9$ and $k/n \to 0$ as $n\to \9$.
Set
 \begin{align} \label{def-Yn*}
    \wt Y_n : = \frac{3\sqrt{2} \pi (\wt{a}(t))^\frac 23}{2} \frac{x_{n-k}-[1-(\frac{k}{\wt{a}(t)n})^{\frac{2}{3}}]}
    {\frac{\sqrt{\log k}}{n^{\frac{2}{3}}k^{\frac{1}{3}}}},
 \end{align}
where $t$ is the unique real number
such that
$t = 1- (\frac{k}{\wt{a}(t) n})^\frac 23
+ \frac{\sqrt{2}}{3\pi (\wt{a}(t))^\frac 23} \frac{\sqrt{\log k}}{n^\frac 23 k^\frac 13} \g_n \xi$.

Then,
for any sequence $\{\g_n\}$ such that $1\ll \g_n \ll \sqrt{\log k}$,
$\{\g^{-1}_n Y_n\}$ satisfies the $LDP(\g_n^2, x^2/2)$.
\end{theorem}

\begin{remark}
The existence and uniqueness of $t$ in Theorem \ref{MDP-Thm} above
can be proved by using contraction mapping arguments on $[1-\delta, 1]$ with $\delta>0$ small enough.
When $k \sim n^\g$ and $0<\g\leq \min\{\frac 25, \frac{1}{2m}\}$,
we can replace $\wt a(t)$ in \eqref{def-Yn*} above by the constant $a_1$ as in Theorem \ref{Edge-Thm}.
See also Remark \ref{Rem-Exp-edge} below.
\end{remark}

\begin{remark}
Theorem \ref{MDP-Thm} is  motivated by the works \cite{DE13,DE13.2},
where the moderate deviation principle of eigenvalues of
Wigner matrices (including the GUE) was proved in the bulk and edge cases.
\end{remark}

By virtue of the determinantal structure of  unitary invariant ensembles,
the proof of Theorems \ref{Bulk-Thm}, \ref{Edge-Thm} and \ref{MDP-Thm} is mainly based on
the central limit theorems in \cite{CL95, S00, S00.2} and the moderate deviation principle in \cite{DE13},
which in turn rely on the asymptotical estimates of
expectation and variance for the counting statistic for an interval.
Such estimates actually can be derived by  the analysis of
corresponding Christoffel-Darboux kernels.

Unlike in \cite{G05, Z14}, it is technically more involved to obtain
the asymptotics of  Christoffel-Darboux kernels $\mathscr{K}_n(x,y)$
for the unitary invariant ensembles considered here,
mainly due to the complicated formulations of
equilibrium density functions
(see \eqref{rhovn} and \eqref{rhov.1} below).
As a matter of fact,
when deriving the estimates of  expectation,
we have to obtain the asymptotics
of $\mathscr{K}_n(x,x)$ in the whole real line, not
just in the interior of  support $(-1,1)$,
because of the lack of
symmetry $\mathscr{K}_n(x,x)=\mathscr{K}_n(-x,-x)$.
Moreover,
for the estimates of   variance,
the straightforward computations
as in  \cite[Lemma 2.3]{G05} and \cite[Lemma
4.1]{Z14} are no longer applicable here to obtain the asymptotics of kernels $\mathscr{K}_n(x,y)$, $x\not =y$.

The key idea to overcome these difficulties is to reformulate
the Christoffel-Darboux kernels in terms of the solutions of Riemann-Hilbert problems
(see \eqref{k-u} and \eqref{5.K.5} below).

In \cite{DG07}-\cite{DKMVZ99} the steepest descent method, introduced by Deift and Zhou
in \cite{DZ93}, has been developed to obtain the asymptotics of
solutions of Riemann-Hilbert problems and has been applied successfully to prove universality for
a variety of statistical quantities arising in unitary invariant ensembles.

In view of the key identities \eqref{k-u} and \eqref{5.K.5} below,
we employ here the Riemann-Hilbert approach
to obtain the crucial asymptotic estimates of Christoffel-Darboux kernels,
which indeed constitute the main
part of  present work and would be also of independent interests.
Once these estimates obtained,
the Gaussian fluctuations and moderate deviations
can be proved by using similar arguments as in \cite{G05, Z14} and \cite{DE13,DE13.2}, respectively.

We would also like to mention that,
the growing variance statistics
may be also derived by adjusting the arguments in \cite{BD14}.

The remainder of this article is organized as follows.
Section \ref{RHM} -- \ref{MDP} are devoted to
unitary invariant ensembles with  Freud-type potentials.
First, in Section
\ref{RHM} we briefly review the Riemann-Hilbert approach developed
in \cite{DKMVZ99},
and then in Section \ref{K-K} we prove the
key asymptotic estimates
of  Christoffel-Darboux kernels. Section \ref{Gauss-Fluct} mainly contains
the proof of  Gaussian fluctuations in Theorems \ref{Bulk-Thm} and \ref{Edge-Thm}.
The precise asymptotics of  expectations and variances are also  given.
Section \ref{MDP} includes the proof of moderate deviations in Theorem \ref{MDP-Thm}.
In Section \ref{Sec-UC} we treat  unitary invariant ensembles with  uniform convex potentials.
For simplicity of
exposition, some technical details
are postponed to the Appendix.

{\it \bf Notation.} Throughout this article, $\#I$ denotes the number of eigenvalues in
the interval $I\subseteq \mathbb{R}$. For two sequence of real numbers $f_n$ and $g_n$, $n\geq 1$, $f_n=\mathcal{O}(g_n)$ means that
$|f_n/g_n|$ stays bounded, and $f_n\ll g_n$ means $\lim_n f_n/g_n =0$.
The notations $C$ and $c$ denote constants which may
change from one line to another.

\section{Riemann-Hilbert approach}\label{RHM}

We start with the Freud-type potential \eqref{v.2}.
Let $Q(x)=\sum_{j=0}^{2m}q_jx^j$,
$q_{2m}=\frac{\Gamma(m)\Gamma(\frac{1}{2})}{\Gamma(\frac{2m+1}{2})}$,
$m\geq 1$. Define the n-th Mhasker-Rakhmanov-Saff numbers
$\alpha_n,\beta_n$ by
\begin{align}
    \frac{1}{2\pi}\int_{\alpha_n}^{\beta_n}\frac{Q'(t)(t-\alpha_n)}{\sqrt{(\beta_n-t)(t-\alpha_n)}}dt=n,  \ \
    \frac{1}{2\pi}\int_{\alpha_n}^{\beta_n}\frac{Q'(t)(\beta_n-t)}{\sqrt{(\beta_n-t)(t-\alpha_n)}}dt =-n.
    \label{mrs2}
\end{align}
It follows from \cite[Proposition $5.2$]{DKMVZ99} that $\alpha_n$ and $\beta_n$ exist for $n$ large enough
and can be expressed in a power series in $n^{-\frac{1}{2m}}$.
Set
\begin{equation*}
   (V(x):=) V_n(x)=\frac{1}{n}Q(c_nx+d_n) \label{vn}
\end{equation*}
with $ c_n =\frac 12 (\beta_n-\alpha_n)$,
$d_n =\frac 12 (\beta_n+\alpha_n)$.
We have  that (\cite[ $(5.17)$, $(5.18)$]{DKMVZ99}),
$ V_n=\sum\limits_{k=0}^{2m}v_{n,k}x^k\ \in\mathbb{P}_{2m}^+,$
where
\begin{align} \label{vnk}
v_{n,2m}=\frac{1}{mA_m}+\mathcal{O}(n^{-\frac{1}{m}}),\ \
v_{n,k}=\mathcal{O}(n^{\frac{k}{2m}-1}),\ \  0\leq k\leq 2m-1,
\end{align}
and $A_m=\prod_{j=1}^m\frac{2j-1}{2j}$.

\begin{remark}
When $Q$ is the monomial polynomial of even degree as in \cite{Z14},
we have that $\beta_n=-\alpha_n=n^\frac{1}{2m}$. Hence,
$c_n=n^\frac{1}{2m}$, $d_n=0$ and  $V_n\equiv Q$.
\end{remark}

We have the following formula for the equilibrium density function.
\begin{theorem} (\cite[Proposition 5.3]{DKMVZ99}) \label{4.rhon}
There exists $N>0$, such that for all $n\geq N$,
\begin{equation}
    \rho_{V_n}(x)=\frac{1}{2\pi}\sqrt{1-x^2}h_n(x) \chi_{[-1,1]}(x) ,\label{rhovn}
\end{equation}
where
\begin{equation}\label{hn}
    h_n(x)=\sum_{k=0}^{2m-2}h_{n,k}x^k,\ \ h_{n,k}=\sum_{j=0}^{[\frac{2m-2-k}{2}]}A_j(k+2+2j)v_{n,k+2+2j}.
\end{equation}
 Furthermore, for some $h_0>0$, $h_n(x)>h_0$
for all $n\geq N$ and $x\in\mathbb{R}$.
\end{theorem}

\begin{lemma} \label{4.F} Let $\delta\in(0,1)$ and $N$ be
as in Theorem \ref{4.rhon}. Then $1/\rho_{V_n}$ and $|\rho'_{V_n}|$
are uniformly bounded for all $n\geq N$ and
$x\in[-1+\delta,1-\delta]$.
\end{lemma}
(See the Appendix for the proof.)

Below we assume that $n$ is large enough such that Theorem
\ref{4.rhon} holds.
Set
\begin{equation}   \label{Fn1}
F_n(x) :=\left|\int_x^1\frac{1}{2\pi}\sqrt{|1-y^2|}h_n(y)dy\right|,
\widetilde{F}_n(x) :=\left|\int_{-1}^x\frac{1}{2\pi}\sqrt{|1-y^2|}h_n(y)dy\right|.
\end{equation}
Then, $F_n(x)=\int_x^1\rho_{V_n}(y)dy$ , $x\in(-1,1)$,
and
$F_n(x)=\widetilde{F}_n(-x)$ if $Q(x)=Q(-x)$.

The $j$-th orthogonal polynomials $p_j(x)$ and the
Christoffel-Darboux kernels $K_j(x,y)$ with respect to the weight
$e^{-Q(x)}$ are defined by
\begin{align*}
    & p_j(x)=\gamma_jx^j+\dots,~~\gamma_j> 0,\ \ j\geq 0, \\
    & \int p_i(x)p_j(x)e^{-Q(x)}dx=\delta_{ij},\ \  i,j\geq 0, \\
    & K_j(x,y)=\sum_{i=0}^{j-1}p_i(x)p_i(y)e^{-\frac{Q(x)+Q(y)}{2}},\ \ j\geq 1.
\end{align*}
For the scaled weight $e^{-nV_n(x)}(=e^{-Q(c_nx+d_n)})$,
we define the $j$-th orthogonal polynomials $p_j(x;n)$ and the corresponding kernels $\mathscr{K}_j(x,y)$
similarly  as follows
\begin{align}
    & p_j(x;n)=\gamma_j^{(n)}\pi_j(x;n), \label{p-gpi}  \\
    & \int p_i(x;n)p_j(x;n)e^{-nV_n(x)}dx=\delta_{ij},\ \ i,j\geq 0, \label{pipj-op} \\
    & \mathscr{K}_j(x,y)=\sum_{i=0}^{j-1}p_i(x;n)p_i(y;n)e^{-n\frac{V_n(x)+V_n(y)}{2}},\ \ j\geq 1,  \label{kj-p}
\end{align}
where  $\gamma_j^{(n)}>0$ and $\pi_j(x;n)$ are monic polynomials.
It is straightforward to verify
\begin{align}
       p_i(x;n)=&\sqrt{c_n}p_i(c_nx+d_n), \label{pi-pi}\\
       \gamma_i^{(n)}=&c_n^{i+\frac{1}{2}}\gamma_i,\ \ i\geq 0, \label{gi-gi}\\
       \mathscr{K}_n(x,y)=&c_nK_n(c_nx+d_n,c_ny+d_n). \label{kn-kn}
\end{align}

Below we recall the Riemann-Hilbert problem and
the steepest descent method, which was introduced by Deift and Zhou
in \cite{DZ93} and later developed in
\cite{DKMVZ99.0}-\cite{DVZ97} to analyze the asymptotics of the
solutions of Riemann-Hilbert problem.

Let $U:\mathbb{C}/\mathbb{R}\rightarrow \mathbb{C}^{2\times2}$ be an analytic
matrix-valued function, which solves the Riemann-Hilbert problem,
\begin{eqnarray*}
    U_+(s)=U_-(s)\left(\begin{array}{ccc}
1&e^{-nV_n(s)}\\{ }&1\end{array}\right), \ s\in\mathbb{R},
\end{eqnarray*}
\begin{equation*}
    U(z)\left(\begin{array}{cc} z^{-n}&{ }\\{
    }&z^n\end{array}\right)=I+\mathcal{O}(\frac{1}{|z|}),\ as\ |z|\rightarrow\infty.
\end{equation*}

The fundamental relation between the solutions of Riemann-Hilbert problem and the
orthogonal polynomials, observed by Fokas, Its and Kitaev
\cite{FIK91}, is that
\begin{equation}
    U_{11}(z)=\frac{1}{\gamma_n^{(n)}}p_n(z;n),\ \ U_{21}(z)=-2\pi
i \gamma_{n-1}^{(n)}p_{n-1}(z;n). \label{up}
\end{equation}

Set
\begin{equation*}
    g_n(z) : =\int_{-1}^1\psi_n(t)\log
(z-t)dt,\ \ z\in\mathbb{C}/(-\infty,1],
\end{equation*}
where
\begin{equation} \label{psi}
    \psi_n(z)=\frac{1}{2\pi}(1-z)^{\frac{1}{2}}(1+z)^{\frac{1}{2}}h_n(z),\ \ z\in\mathbb{C}/((-\infty,-1]\cup
[1,\infty))
\end{equation}
with the analytic branch chosen by $\arg(1-x)=\arg(1+x)=0$,
$x\in(-1,1)$. Let
\begin{equation*}
    \xi_n(z) : =-2\pi
i\int_1^z\psi_n(y)dy,\ \ z\in\mathbb{C}/(-\infty,-1]\cup [1,\infty).
\end{equation*}
We have that (\cite[(8.29)]{DKMVZ99})
\begin{equation} \label{gn}
    g_n(z)=\frac{1}{2}(V_n(z)+l_n+\xi_n(z)), \ \ z\in\mathbb{C}^{+},
\end{equation}
where $l_n$ is same as $l$ in \eqref{rho-vl1} with $V_n$ replacing $V$ there.

Using the Pauli matrix $\sigma_3=\left(
                                   \begin{array}{cc}
                                     1 & 0 \\
                                     0 & -1 \\
                                   \end{array}
                                 \right)
$, we set
\begin{equation}
    T(z) : =e^{-n\frac{l_n}{2}\sigma_3}U(z)e^{-n(g_n(z)-\frac{l_n}{2})\sigma_3},
    \ z\in \mathbb{C}/\mathbb{R}.
    \label{T}
\end{equation}
and
\begin{equation} \label{S}
    S(z): =\left\{
        \begin{array}{ll}
          T(z), & \hbox{outside the lens-shaped region;} \\
          T(z)\left(\begin{array}{cc} 1&0\\-e^{-n\xi_n}&1\end{array}\right), & \hbox{in the upper lens region;} \\
          T(z)\left(\begin{array}{cc} 1&0\\e^{n\xi_n}&1\end{array}\right), & \hbox{in the lower lens region,}
        \end{array}
      \right.
\end{equation}
with the lens regions  as in \cite[fig. $6.1$]{DKMVZ99}.

Next, we recall the delicate paramatrices $P_n$ in the small balls $U_{\pm
1}$ centered on $\pm 1$ with
the radius $\delta$ sufficiently small, respectively. Define
the functions $f_n$ and $\widetilde{f}_n$ in $U_{1}$ and
$U_{-1}$ respectively by
\begin{align}
    & (-f_n(z))^{\frac{3}{2}}=-n\frac{3\pi}{2}\int_1^z\psi_n(y)dy,\ \
z\in U_1/[1,\infty), \label{fn1} \\
    & (\widetilde{f}_n(z))^{\frac{3}{2}}=n\frac{3\pi}{2}\int_{-1}^z\psi_n(y)dy,\ \
z\in U_{-1}/(-\infty,-1]. \label{f^n1}
\end{align}
We have that (see $(7.14)$, $(7.21)$, $(7.38)$, $(7.36)$ and
$(7.37)$ in \cite{DKMVZ99}),
\begin{align}
    \frac{2}{3}(f_n(z))^{\frac{3}{2}}=n\varphi_n(z),\ or,\ f_n(z)=n^{\frac{2}{3}}(z-1)(\widehat{\phi}_n(z))^{\frac{2}{3}}, \label{fn2} \\
    \frac{2}{3}(-\widetilde{f}_n(z))^{\frac{3}{2}}=n\widetilde{\varphi}_n(z),\ or,\ \widetilde{f}_n(z)=n^{\frac{2}{3}}(z+1)(\widehat{\widetilde{\phi}}_n(z))^{\frac{2}{3}}, \label{f^n2}
\end{align}
where
\begin{equation}\label{phin}
 \varphi_n(z)=\left\{
       \begin{array}{ll}
         -\frac{1}{2}\xi_n(z)=\pi i \int_1^z\psi_n(y)dy, & \hbox{$z\in\mathbb{C}^+$;} \\
         \frac{1}{2}\xi_n(z)=-\pi i\int_1^z\psi_n(y)dy , & \hbox{$z\in\mathbb{C}^-$;}
       \end{array}
     \right.
\end{equation}
\begin{equation}\label{phi^n}
 \widetilde{\varphi}_n(z)=\left\{
       \begin{array}{ll}
         \varphi_n(z)+\pi i=\pi i \int_{-1}^z\psi_n(y)dy, & \hbox{$z\in\mathbb{C}^+$;} \\
         \varphi_n(z)-\pi i=-\pi i\int_{-1}^z\psi_n(y)dy , & \hbox{$z\in\mathbb{C}^-$;}
       \end{array}
     \right.
\end{equation}
and $\widehat{\phi}_n$, $\widehat{\widetilde{\phi}}_n$ are analytic
functions in $U_{1}$ and $U_{-1}$, respectively.

The paramatrices $P_n$ in $U_{\pm 1}$ are defined as follows.
\\$(i)$. In the region $U_1/f_n^{-1}(\gamma_{\sigma})$ with the contour $\gamma_{\sigma}$
as in \cite[fig. $7.1$]{DKMVZ99}, set

\begin{equation}
    P_n:=E_n\Psi^{\sigma}(f_n)e^{n\varphi_n\sigma_3}, \label{Pn1}
\end{equation}
where $E_n=\sqrt{\pi}e^{\frac{\pi i}{6}}\left(
                                          \begin{array}{cc}
                                            1 & -1 \\
                                            -i & -i \\
                                          \end{array}
                                        \right)
\left(
  \begin{array}{cc}
    H_n &   \\
      & H_n^{-1} \\
  \end{array}
\right) $, $H_n=f_n^{\frac{1}{4}}a^{-1}$, and
\begin{equation} \label{Psi}
\Psi^{\sigma}(z)=\left\{
  \begin{array}{ll}
    AI(z)e^{-\frac{\pi i}{6}\sigma_3}, & \hbox{$z\in I: 0<\arg z<\frac{2\pi}{3}$;} \\
    AI(z)e^{-\frac{\pi i}{6}\sigma_3}\left(
                                     \begin{array}{cc}
                                       1 & 0 \\
                                       -1 & 1 \\
                                     \end{array}
                                   \right)
, & \hbox{$z\in II: \frac{2\pi}{3}< \arg z<\pi$;} \\
    \widetilde{AI}(z)e^{-\frac{\pi i}{6}\sigma_3}\left(
                                     \begin{array}{cc}
                                       1 & 0 \\
                                       1 & 1 \\
                                     \end{array}
                                   \right), & \hbox{$z\in III: -\pi < \arg z <-\frac{2\pi}{3}$;} \\
    \widetilde{AI}(z)e^{-\frac{\pi i}{6}\sigma_3}, & \hbox{$z\in IV: -\frac{2\pi}{3}<\arg z<0$.}
  \end{array}
\right.
\end{equation}
Here, $AI(z)$ and $\widetilde{AI}(z)$ denote $\left(
                                    \begin{array}{cc}
                                      Ai(z) & Ai(\omega^2 z) \\
                                      Ai'(z) & \omega^2Ai'(\omega^2 z)
                                    \end{array}
                                  \right)
$, $\left(
                                    \begin{array}{cc}
                                      Ai(z) & -\omega^2 Ai(\omega z) \\
                                      Ai'(z) & -Ai'(\omega z)
                                    \end{array}
                                  \right)$,
respectively, $\omega=e^{\frac{2\pi
                                  i}{3}}$, and
$Ai$ is the Airy function, uniquely determined by
the equation $Ai''(z)=z Ai(z)$ with $\lim_{x\to \infty}
\sqrt{4\pi}
x^{\frac{1}{4}} e^{\frac{2}{3}x^{\frac{3}{2}}} Ai(x)=1$.

$(ii)$. In the region
$U_{-1}/\widetilde{f}_n^{-1}(\widetilde{\gamma}_{\sigma})$ with the
contour $\widetilde{\gamma}_{\sigma}$  as in \cite[fig.
$7.3$]{DKMVZ99}, set
\begin{equation} \label{Pn-1}
    P_n:=\widetilde{E}_n\widetilde{\Psi}^{\sigma}(\widetilde{f}_n)e^{n\widetilde{\varphi}_n\sigma_3}
\end{equation}
with
$\widetilde{\Psi}^{\sigma}(z)=\sigma_3\Psi^{\sigma}(-z)\sigma_3$,
$\widetilde{E}_n=\sqrt{\pi}e^{\frac{\pi i}{6}}\left(
                                                    \begin{array}{cc}
                                                      1 & 1 \\
                                                      i & -i \\
                                                    \end{array}
                                                  \right)
\left(
  \begin{array}{cc}
    \widetilde{H}_n &   \\
      & \widetilde{H}_n^{-1} \\
  \end{array}
\right) $, $\widetilde{H}_n=(-\widetilde{f}_n)^{\frac{1}{4}}a$.
 \\

Finally, set
\begin{equation} \label{R}
    R:=\left\{
        \begin{array}{ll}
          SP_n^{-1}, & \hbox{for $z\in U_1\cup U_{-1}$;} \\
          SN^{-1}, & \hbox{otherwise,}
        \end{array}
      \right.
\end{equation}
where
\begin{equation} \label{N}
N=\frac{1}{2}\left(
                       \begin{array}{cc}
                         a+a^{-1} & i(a^{-1}-a) \\
                         i(a-a^{-1}) & a+a^{-1} \\
                       \end{array}
                     \right),
\end{equation}
and
\begin{equation} \label{a}
       a(z)=(\frac{z-1}{z+1})^{\frac{1}{4}},\ z\in\mathbb{C}/[-1,1]
\end{equation}
with the analytic branch chosen by $\arg(x-1)=\arg(x+1)=0$, for
$x>1$.  We have the  asymptotic
expansions of $R$ below (see \cite[(7.64)]{DKMVZ99}, \cite [(3.6), (3.7)]{DG07}),
\begin{align}  \label{4.R}
    R(z) = I+\frac{1}{n}\sum_{k=0}^{\infty}r_k(z)
    n^{-\frac{k}{2m}}, \ \
    \frac{d}{dz}R(z) = \frac 1 n \sum\limits_{k=0}^{\infty}\frac{d}{dz}r_k(z)n^{-\frac{k}{2m}},
\end{align}
where $r_k(z)$, $\frac{d}{dz}r_k(z)$, $0\leq k<\9$, are bounded
functions and analytic in the complement of set $\partial U_1
\cup \partial U_{-1}$, and these expansions are uniform for
$z\in\mathbb{C}/\widehat{\Sigma}_R$ with $\widehat{\Sigma}_R$ as
in \cite[fig.$7.6$]{DKMVZ99}.

\begin{remark}
If $z=x \in \mathbb{R}$, we take the limiting expressions as
$z$ is approaching from the upper half-plane. Thus, if $x>1$,
$\psi_n(x)$ means $\lim_{\epsilon\rightarrow
0^+}\psi_n(x+i\epsilon)$.
\end{remark}

\section{Asymptotics of Christoffel-Darboux kernels} \label{K-K}

This section is mainly devoted to the asymptotics  of
Christoffel-Darboux kernels corresponding to
Freud-type potentials.

\begin{lemma} \label{4.K.0}
Take any sufficiently small $\delta>0$, we have

$(i)$. For $x\in (-1+\delta,1-\delta)$,
\begin{equation}\label{4.K.01}
    \mathscr{K}_n(x,x)=n\rho_{V_n}(x)+\mathcal {O}(1).
\end{equation}

$(ii)$. For  $x\in(1-\delta,1+\delta)$,
\begin{align}\label{4.K.02}
    \mathscr{K}_n(x,x)
    =&\left[\frac{1}{4}\frac{f_n'(x)}{f_n(x)}-\frac{a'(x)}{a(x)}\right]2Ai(f_n(x))Ai'(f_n(x))  \nonumber  \\
    &+f_n'(x)\left[(Ai')^2(f_n(x))-f_n(x)Ai^2(f_n(x))\right]
    +\mathcal {O}(n^{-\frac{5}{6}}).
\end{align}

$(iii)$. For  $x\in(-1-\delta,-1+\delta)$,
\begin{align}\label{4.K.03}
    \mathscr{K}_n(x,x)=&-\left[\frac{1}{4}\frac{\widetilde{f}'_n(x)}{\widetilde{f}_n(x)}+\frac{a'(x)}{a(x)}\right]2Ai(-\widetilde{f}_n(x))Ai'(-\widetilde{f}_n(x))\nonumber \\
            &+\widetilde{f}'_n(x)\left[(Ai')^2(-\widetilde{f}_n(x))+\widetilde{f}_n(x)Ai^2(-\widetilde{f}_n(x))\right]+\mathcal
            {O}(n^{-\frac{5}{6}}).
\end{align}

$(iv)$. For  $x\in \mathbb{R}/(-1-\delta,1+\delta)$,
\begin{equation}\label{4.K.04}
    \mathscr{K}_n(x,x)
    =\frac{1}{4\pi}\frac{1}{(x-1)(x+1)}e^{-2n\varphi_n(x)}+\mathcal {O}(n^{-1}).
\end{equation}
\end{lemma}

{\it \bf Proof.} First  by  the Christoffel-Darboux formula
(cf. \cite[$(3.48)$]{D99}) and  (\ref{up}),
\begin{equation}\label{k-u}
    2\pi
    i(x-y)\mathscr{K}_n(x,y)=(1,0)U(x)^TU(y)^{-T}(0,1)^Te^{-n\frac{V_n(x)+V_n(y)}{2}}.
\end{equation}
This identity is the key to relate the Christoffel-Darboux  kernels with the solutions of Riemann-Hilbert problem
and so
enables us to employ the
Riemann-Hilbert approach to obtain the asymptotics  of these kernels.

$(i)$. For $x,y\in (-1+\delta,1-\delta)$, by (\ref{T}) and
(\ref{S}),
\begin{align} \label{k-u*}
    U=e^{\frac{n}{2}l_n\sigma_3}S\left(
                                    \begin{array}{cc}
                                      1 & 0 \\
                                      e^{-n\xi_n} & 1 \\
                                    \end{array}
                                  \right)
    e^{n(g_n-\frac{l_n}{2})\sigma_3}.
\end{align}
Then, by (\ref{k-u}), (\ref{gn}) and (\ref{phin}), direct
computations show that
\begin{align} \label{kn+error1***}
    2\pi i(x-y)\mathscr{K}_n(x,y)
    =&(e^{-n\varphi_n(x)},e^{n\varphi_n(x)})S(x)^TS(y)^{-T}(-e^{n\varphi_n(y)},e^{-n\varphi_n(y)})^T.
\end{align}

In order to obtain the leading term of the right-hand side above, we
note that $S^T(x)=S^T(y)+(x-y) \Delta_S(x,y)$, where \begin{align}
\label{kn+error1*}
 \Delta_{S}(x,y)=\int_0^1(S^T)'(y+t(x-y))dt.
 \end{align}
This yields that
\begin{equation} \label{kn+error1**}
    S^T(x)S^{-T}(y)=Id+(x-y) \Delta_S(x,y) S^{-T}(y).
\end{equation}
Then, plugging \eqref{kn+error1**} into \eqref{kn+error1***},
since $\varphi_n(x)=-\pi i F_n(x)$, $x\in(-1,1)$, we obtain
\begin{equation}\label{kn+error1}
    2\pi i(x-y) \mathscr{K}_n(x,y)=-2 i \sin [n\pi
    (F_n(x)-F_n(y))] +(x-y) I_1(x,y),
\end{equation}
where
$I_1(x,y):=(e^{-n\varphi_n(x)},e^{n\varphi_n(x)})[\Delta_{S}(x,y)S^{-T}(y)](-e^{n\varphi_n(y)},e^{-n\varphi_n(y)})^T.$

Thus, taking $y=x$ we get
\begin{align*}
    &2\pi i\mathscr{K}_n(x,x) \nonumber \\
    =&2\pi i n \rho_{V_n}(x)
    +(e^{-n\varphi_n(x)},e^{n\varphi_n(x)})\left[(S^T)'(x)S^{-T}(x)\right](-e^{n\varphi_n(x)},e^{-n\varphi_n(x)})^T.
\end{align*}
In view of \eqref{R} -- \eqref{4.R}, $S(x)$
and $S'(x)$ are uniformly bounded for $x\in[-1+\delta,1-\delta]$,
hence \eqref{4.K.01} follows.

$(ii)$. For $x,y\in (1-\delta,1),\ or,\ x,y\in (1,1+\delta)$,
similar calculations show that
\begin{align} \label{kn-en}
    2\pi i(x-y)\mathscr{K}_n(x,y)
    =&e^{-\frac{\pi
i}{3}}(1,0)[AI(f_n(x))]^TE_n^T(x)R^T(x) \nonumber \\
     &\quad \cdot R^{-T}(y)E_n^{-T}(y)[AI(f_n(y))]^{-T}(0,1)^T.
\end{align}
(see the Appendix for the proof.)

Regarding the leading term of the right-hand side above, using
\eqref{kn+error1**} with $S$ replaced by $R$, we obtain
\begin{align} \label{kn-error2}
     2\pi i(x-y)\mathscr{K}_n(x,y)
    =&e^{-\frac{\pi i}{3}}(1,0)[AI(f_n(x))]^TE_n^T(x)
       E_n^{-T}(y)[AI(f_n(y))]^{-T}(0,1)^T \nonumber  \\
    &\ +(x-y)e^{-\frac{\pi
    i}{3}}I_2(x,y),
\end{align}
where
\begin{align}   \label{I2}
    I_2(x,y):=&(1,0)[AI(f_n(x))]^TE_n^T(x) \nonumber \\
    &\quad \cdot  \Delta_{R}(x,y)R^{-T}(y)E_n^{-T}(y)[AI(f_n(y))]^{-T}(0,1)^T.
\end{align}

Then, using the expressions of $AI$, $E_n$ and the asymptotics
\eqref{4.R},  we have that
$I_2(x,y)$ is of order $n^{-\frac{5}{6}}$ and
\begin{align} \label{kn-ai}
    2\pi i(x-y)\mathscr{K}_n(x,y) &
   =  (-2\pi
   i)\bigg[-Ai(f_n(x))Ai'(f_n(y))\frac{f_n^{\frac{1}{4}}(x)}{f_n^{\frac{1}{4}}(y)}\frac{a(y)}{a(x)}
   \nonumber \\
   &
   +Ai'(f_n(x))Ai(f_n(y))\frac{f_n^{\frac{1}{4}}(y)}{f_n^{\frac{1}{4}}(x)}\frac{a(x)}{a(y)}\bigg]
   +(x-y)\mathcal {O}(n^{-\frac{5}{6}}).
\end{align}
The proof is postponed to the Appendix.

Hence, taking the Taylor expansion and using
$Ai''(x)=xAi(x)$  we obtain
\begin{align*}
   & 2\pi i(x-y)\mathscr{K}_n(x,y)\\
   =&-2\pi
   i(y-x)\bigg\{\left[\frac{1}{4}\frac{f_n'(x)}{f_n(x)}-\frac{a'(x)}{a(x)}\right]2Ai(f_n(x))Ai'(f_n(x))\\
   &-f_n(x)f_n'(x)Ai^2(f_n(x))+f_n'(x)(Ai')^2(f_n(x))\bigg\}+\mathcal{O}((y-x)^2)+(x-y)\mathcal{O}(n^{-\frac{5}{6}}),
\end{align*}
which implies  (\ref{4.K.02}).

$(iii)$. For $x,y\in(-1-\delta,-1),\ or,\ x,y \in (-1,-1+\delta)$,
the proofs are similar to those in the previous case. First we compute that
\begin{align} \label{kn-en2}
    2\pi i(x-y)\mathscr{K}_n(x,y)
    =& (-1)e^{-\frac{\pi
    i}{3}}(1,0)[\widetilde{AI}(-\widetilde{f}_n(x))]^{T}\sigma_3\widetilde{E}_n^{T}(x)R^T(x) \nonumber \\
    &\quad \cdot R^{-T}(y)\widetilde{E}_n^{-T}(y)\sigma_3^{-1}[\widetilde{AI}(-\widetilde{f}_n(y))]^{-T}(0,1)^T.
\end{align}
(See the  Appendix for the proof.)

Then, using \eqref{kn+error1**} with $S$ replaced by $R$ we have
\begin{align} \label{4.K.(3).2}
    &2\pi i(x-y)\mathscr{K}_n(x,y) \nonumber  \\
    =&(-1)e^{-\frac{\pi i}{3}}(1,0)[\widetilde{AI}(-\widetilde{f_n}(x))]^{T}\sigma_3\widetilde{E}_n^{T}(x)
       \widetilde{E}_n^{-T}(y)\sigma_3^{-1}[\widetilde{AI}(-\widetilde{f}_n(y))]^{-T}(0,1)^T \nonumber \\
    &\quad  -e^{-\frac{\pi i}{3}}(x-y)I_3(x,y),
\end{align}
where
\begin{align} \label{I3}
 I_3(x,y)
:=&(1,0)[\widetilde{AI}(-\widetilde{f}_n(x))]^{T}\sigma_3\widetilde{E}_n^{T}(x)\nonumber
\\
  &\quad \cdot  \Delta_{R}(x,y)R^{-T}(y)\widetilde{E}_n^{-T}(y)\sigma_3^{-1}[\widetilde{AI}(-\widetilde{f}_n(y))]^{-T}(0,1)^T.
\end{align}
Using the expressions of $\wt{AI}$, $\wt f_n$ and arguing as
in the case $(ii)$  we get that
    \begin{align}\label{kn-error3esti}
        &2\pi i (x-y) \mathscr{K}_n(x,y) \nonumber  \\
        =&(-2\pi i)\bigg[Ai(-\widetilde{f}_n(x))Ai'(-\widetilde{f}_n(y))\frac{\widetilde{f}_n(x)^{\frac{1}{4}}}{\widetilde{f}_n(y)^{\frac{1}{4}}}\frac{a(x)}{a(y)}
       \nonumber \\
        &\qquad  \qquad -Ai'(-\widetilde{f}_n(x))Ai(-\widetilde{f}_n(y))\frac{\widetilde{f}_n(y)^{\frac{1}{4}}}{\widetilde{f}_n(x)^{\frac{1}{4}}}\frac{a(y)}{a(x)}\bigg]
       + (x-y) \mathcal{O}(n^{-\frac{5}{6}}).
    \end{align}

Consequently, the Taylor expansion yields that
\begin{align*}
    &2\pi i(x-y)\mathscr{K}_n(x,y)\\
    =&-2\pi i(y-x) \bigg\{-\left[\frac{1}{4}\widetilde{f}_n^{-1}(x)\widetilde{f}'_n(x)+\frac{a'(x)}{a(x)}\right]2Ai(-\widetilde{f}_n(x))Ai'(-\widetilde{f}_n(x))\\
    &+\widetilde{f}'_n(x)\left[\widetilde{f}_n(x)Ai^2(-\widetilde{f}_n(x))+(Ai')^2(-\widetilde{f}_n(x))\right]\bigg\}
    +\mathcal{O}(y-x)^2 +(x-y)\mathcal{O}(n^{-\frac{5}{6}}),
\end{align*}
which implies (\ref{4.K.03}).

$(iv)$. For $x,y\in \mathbb{R}/(-1-\delta,1+\delta)$, by
\eqref{T} and \eqref{S},
\begin{equation*}
     U=e^{n\frac{l_n}{2}\sigma_3}Se^{n(g_n-\frac{l_n}{2})\sigma_3},
\end{equation*}
which along with \eqref{k-u} and \eqref{gn}  implies  that
$$ 2\pi i(x-y)\mathscr{K}_n(x,y)=e^{-n(\varphi_n(x)+\varphi_n(y))}(1,0)S^T(x)S^{-T}(y)(0,1)^T.$$

Then,
using $S=RN$ and applying \eqref{kn+error1**} twice with $S$
replaced by $R$ and $N$, respectively, we obtain
\begin{align} \label{kn+error4}
    2\pi i(x-y)\mathscr{K}_n(x,y)
    =&(x-y)e^{-n(\varphi_n(x)+\varphi_n(y))}(1,0)
      \Delta_N(x,y)N^{-T}(y)(0,1)^T \nonumber \\
     &+(x-y)I_4(x,y),
\end{align}
where
$$I_4(x,y):=e^{-n(\varphi_n(x)+\varphi_n(y))}(1,0)N^T(x)\Delta_R(x,y)R^{-T}(y)N^{-T}(y)(0,1)^T .$$
Taking into account \eqref{4.R} and that $N(x)$ and $e^{-n\varphi_n(x)}$ are bounded for
$x\in\mathbb{R}/(-1-\delta,1+\delta)$,  we obtain that $I_4(x,y)=\mathcal
{O}(n^{-1})$.
Hence,
\begin{eqnarray*}
    \mathscr{K}_n(x,y)=\frac{1}{2\pi
    i}e^{-n(\varphi_n(x)+\varphi_n(y))}(1,0)\Delta_N(x,y)N^{-T}(y)(0,1)^T+\mathcal
    {O}(n^{-1}),
\end{eqnarray*}
and
\begin{equation*}
    \mathscr{K}_n(x,x)=\frac{1}{2\pi
    i}e^{-2n\varphi_n(x)}(1,0)(N^T)'(x)N^{-T}(x)(0,1)^T+\mathcal
    {O}(n^{-1}),
\end{equation*}
which consequently yields  (\ref{4.K.04}) by   \eqref{N} and \eqref{a}.
The proof  is  complete. \hfill $\square$

Lemmas \ref{4.K.1} and \ref{4.K.1*} below are concerned with  the asymptotics of
kernels $\mathscr{K}_n(x,y)$, $x\not =y$,
in the bulk and edge cases, respectively.

\begin{lemma} \label{4.K.1}
$(i)$. Let $t\in(-1,1)$  and set $\Gamma_1^1 :=\{(x,y):t\leq x\leq
t+\frac{1-t}{\log n},t-\frac{1+t}{\log n}\leq y\leq
t-\frac{1}{n}\}$. Then, for $(x,y)\in\Gamma_1^1$,
\begin{equation}
    \mathscr{K}_n(x,y)=\frac{\sin[\pi n (F_n(y)-F_n(x))]+\mathcal{O}(\frac{1}{\log
n})}{\pi(x-y)},\label{4k6}
\end{equation}
where $F_n(x)$ is defined as in (\ref{Fn1}).

$(ii)$. Let $\delta>0$ be sufficiently small.
Then, for $x,y
\in[-1+\delta,1-\delta]$,
\begin{equation}
    \mathscr{K}^2_n(x,y)=\mathcal{O}(\frac{1}{(x-y)^2}). \label{k7}
\end{equation}

$(iii)$. Let $t\in(-1,1)$.
Then, for $(x,y)\in\{(x,y):t\leq x\leq
t+\frac{1}{n},t-\frac{1}{n}\leq y\leq t\}$,
\begin{equation} \label{k7*}
    \mathscr{K}_n(x,y)=\mathcal{O}(n).
\end{equation}
\end{lemma}

{\it \bf Proof.} By (\ref{kn+error1}), we have
\begin{equation} \label{kfn}
     2\pi i(x-y)\mathscr{K}_n(x,y)
    =-2i\sin[n\pi(F_n(x)-F_n(y))]+\mathcal{O}(|x-y|),
\end{equation}
which immediately implies \eqref{4k6} and  \eqref{k7}.
As regards $(iii)$, by \eqref{kfn},
\begin{align*}
   \mathscr{K}_n(x,y)
   =& (-1) \frac{\sin[n\pi(F_n(x)-F_n(y))]}{\pi (x-y)} + \mathcal{O}(1) \\
   =& (-1) \frac{\sin[n\pi(F_n(x)-F_n(y))]}{n\pi(F_n(x)-F_n(y))} \frac{n(F_n(x)-F_n(y))}{x-y}   + \mathcal{O}(1).
\end{align*}
Since  $\sup_{x\in \mathbb{R}}|\frac{\sin x}{x}| =\mathcal{O}(1)$,
and by Lemma \ref{4.F},
for some $\xi \in (x,y)$,
\begin{align*}
   \frac{n(F_n(x)-F_n(y))}{x-y}
   =   \frac{nF'_n(\xi)(x-y)}{x-y}
   = - n\rho_{V_n}(\xi)
   =\mathcal{O}(n),
\end{align*}
we obtain \eqref{k7*} and finish the proof.
\hfill $\square$

\begin{lemma} \label{4.K.1*}
For $x
 \in (1-\delta,1+\delta)$, $
 y\in (1-\delta,1)$ with $\delta>0$ sufficient small,
\begin{align} \label{asy-edge-var.1}
    \mathscr{K}_n(x,y)
    =& \frac{1}{x-y}
    \bigg[Ai(f_n(x))Ai'(f_n(y))\frac{f_n^{\frac{1}{4}}(x)}{f_n^{\frac{1}{4}}(y)}\frac{a(y)}{a(x)}  \nonumber \\
    &\qquad \quad -Ai'(f_n(x))Ai(f_n(y))\frac{f_n^{\frac{1}{4}}(y)}{f_n^{\frac{1}{4}}(x)}\frac{a(x)}{a(y)}\bigg]
     + \mathcal{O}(n^{-\frac{5}{6}}).
\end{align}
\end{lemma}

{\bf Proof.} In view of (\ref{kn-ai}), we only need to
prove \eqref{asy-edge-var.1} for  $x\in(1,1+\delta)$
and $y\in(1-\delta,1)$.
For $x\in(1,1+\delta)$, by \eqref{T}, \eqref{S},
 \eqref{Pn1} and \eqref{R},
\begin{align} \label{asy-edge-var.1.1}
    U(x) = e^{n\frac{l_n}{2}\sigma_3} S(x) e^{n(g_n-\frac{l_n}{2})\sigma_3},
\end{align}
where $S(x)=R(x)E_n(x) [AI(f_n(x))]e^{-\frac{\pi i}{6}\sigma_3}
    e^{n\varphi_n\sigma_3}$.
Moreover, for $y\in(1-\delta,1)$,
\begin{align} \label{asy-edge-var.1.2}
  U(y) = e^{n\frac{l_n}{2}\sigma_3}S(y)
        \left(
          \begin{array}{cc}
            1 & 0 \\
            e^{-n\xi_n} & 1 \\
          \end{array}
        \right)
        e^{n(g_n-\frac {l_n}{2})\sigma_3}
\end{align}
with $S(y)=R(y)E_n(y) [AI(f_n(y))] e^{-\frac{\pi i}{6}\sigma_3}
\left(
  \begin{array}{cc}
    1 & 0 \\
    -1 & 1 \\
  \end{array}
\right) e^{n\varphi_n\sigma_3} $.

Then, plugging \eqref{asy-edge-var.1.1} and
\eqref{asy-edge-var.1.2} into \eqref{k-u}  we obtain
\begin{align}
     2\pi i(x-y)\mathscr{K}_n(x,y)
    =&e^{-\frac{\pi
i}{3}}(1,0)[AI(f_n(x))]^TE_n^T(x)R^T(x)
 \nonumber \\
     &\quad \cdot R^{-T}(y)E_n^{-T}(y)[AI(f_n(y))]^{-T}(0,1)^T,
\end{align}
which has the same expressions as in \eqref{kn-en}.

Thus,
similar arguments there yield \eqref{asy-edge-var.1} for
$x\in(1,1+\delta)$, $y\in(1-\delta,1)$.  \hfill $\square$

We conclude this section with the estimates of  orthogonal polynomials.
\begin{lemma} \label{4.B}
There exists a $\delta_0>0$, such that for all
$0<\delta\leq\delta_0$, we have \\
$(i)$. For $x\in \mathbb{R}/(-1-\delta,1+\delta)$,
\begin{equation*}
    p_n(x;n)e^{-\frac{n}{2}V_n(x)}
    \leq C  \bigg[ e^{-n\pi F_n(x)} \chi_{(1+\delta,\9)}(x)
    +  e^{-n\pi\widetilde{F}_n(x)} \chi_{(-\9,-1-\delta)}(x) \bigg].
\end{equation*}
$(ii)$. For $x\in(-1-\delta,1+\delta)$,
\begin{align*}
    p_n(x;n)e^{-\frac{n}{2}V_n(x)} \leq&C\left[1+\frac{1}{|1-x|^{\frac{1}{4}}} \chi_{(1-\delta,1+\delta)}(x)
       +\frac{1}{|1+x|^{\frac{1}{4}}} \chi_{(-1-\delta,-1+\delta)}(x)\right].
\end{align*}
Here $C$ is independent of $n$,
and $\chi_I$ means the characteristic function of $I\subseteq \mathbb{R}$.
\end{lemma}

{\it \bf Proof.} This follows from the Plancherel-Rotach-type
asymptotics of $p_n(x;n)$ in \cite{DKMVZ99} and the asymptotics
of  Airy functions (\cite[(2.60), (2.61), (3.6),
(3.7)]{S06}). \hfill $\square $

\section{Gaussian fluctuations} \label{Gauss-Fluct}

In this section we prove Theorems \ref{Bulk-Thm} and \ref{Edge-Thm}
for Freud-type potentials.

\subsection{Bulk case.} \label{Bulk-Section}

We start with the asymptotics  of the expectation and variance in Propositions \ref{4.exp}
and \ref{VIn-est} below, respectively.

\begin{proposition} \label{4.exp}
Let $t=t(k,n)$ be  as in theorem \ref{Bulk-Thm}. Fix
$\xi\in\mathbb{R}$
and set $a_n :=\frac{\sqrt{\log
n}}{\sqrt{2\pi^2}n\rho_{V_n}(t)},$ $t_n :=t+a_n\xi$ and $
I_n :=[t_n,\infty)$. Then,
\begin{equation} \label{EIn}
    \mathbb{E}_n(\# I_n)=n-k-\frac{\sqrt{\log n}}{\sqrt{2\pi^2}}\xi+\mathcal
    {O}(1).
\end{equation}
\end{proposition}

{\it \bf Proof.} First, since
$ |t(k,n)|= |G^{-1}(k/n)|<1$ for $n$ large enough, Lemma
\ref {4.F} implies that $1/\rho_{V_n}(t)$ are
uniformly bounded and $a_n=\mathcal{O}(\sqrt{\log n} /n)$.

Using the estimates of $\mathscr{K}_n(x,x)$ in Lemma \ref{4.K.0} $(i)$ and  $(iv)$  we have
\begin{align} \label{exp-kn}
    \mathbb{E}_n(\# I_n)
    = \int_{t_n}^{\infty}\mathscr{K}_n(x,x)dx
    = \int_{t_n}^{1-\delta}
    n\rho_{V_n}(x)dx+\int_{1-\delta}^{1+\delta}\mathscr{K}_n(x,x)dx+\mathcal{O}(1).
\end{align}

Moreover, using Lemma \ref{4.K.0}
$(ii)$ and arguing as in the proof of \cite[Lemma 2]{S06} (see also the proof of \cite[(2.2.4)]{Z14.2}), we have
\begin{equation} \label{state}
    \int_{1-\delta}^{1+\delta}\mathscr{K}_n(x,x)dx=\int_{1-\delta}^1n\rho_{V_n}(x)dx+\mathcal{O}(1).
\end{equation}

Therefore, plugging (\ref{state}) into \eqref{exp-kn} and using the Taylor expansion  we
obtain
\begin{align} \label{compute-Exp}
    \mathbb{E}_n(\# I_n)
    =&\int_{t_n}^1n\rho_{V_n}(x)dx+\mathcal{O}(1) \nonumber  \\
    =&n-n\int_{-1}^{t}\rho_{V_n}(x)dx-n\int_{t}^{t+a_n\xi}\rho_{V_n}(x)dx+\mathcal{O}(1) \nonumber  \\
    =&n-k-n\left[\rho_{V_n}(t)a_n\xi + \frac{1}{2}\rho'_{V_n}(\eta)(a_n
    \xi)^2\right]+\mathcal{O}(1)
\end{align}
with $\eta \in (t,t+a_n \xi)$.
Using Lemma \ref{4.F}
we obtain (\ref{EIn})
and finish the proof. \hfill $\square$

\begin{proposition} \label{VIn-est}
 Let $\{t_i\}_{i=1}^\infty$ be a sequence such that $\sup_n|t_n|<1$.
 Set $I_n:=[t_n,\infty),n\in\mathbb{N}$. Then,
\begin{eqnarray} \label{VIn}
    Var_n(\# I_n)=\frac{1}{2\pi^2}\log n
+\mathcal{O}(\log\log n).
\end{eqnarray}
\end{proposition}

{\it \bf Proof.} The arguments are similar to those
in the proof of \cite[Lemma
$3.2$]{Z14}
(see also \cite[Lemma $2.3$]{G05},  \cite[Proposition 2.8]{Z14.2}),
so we give a
sketch of it below.

First, we have
\begin{align} \label{VIn*}
   Var_n(\# I_n)=&\iint\limits_{\Omega_n}\mathscr{K}_n^2(x,y)dxdy \nonumber \\
             =&\iint\limits_{\Gamma}\mathscr{K}_n^2(x,y)dxdy+\iint\limits_{\Omega_n/\Gamma}\mathscr{K}_n^2(x,y)dxdy,
\end{align}
where $\Omega_n=\{(x,y):t_n\leq x< \infty, -\infty<y\leq t_n\}$ and
$\Gamma = \{ (x,y): t_n \leq x \leq 1-\delta, -1+\delta \leq y \leq
t_n \}$.

By virtue of the asymptotic estimates of $\mathscr{K}_n(x,y)$ in Lemma
\ref{4.K.1}, we have
\begin{equation} \label{VIn**}
    \iint\limits_{\Gamma}\mathscr{K}_n^2(x,y)dxdy=\frac{1}{2\pi^2}\log n
+\mathcal{O}(\log\log n),
\end{equation}
where the leading term comes from
the integration on  $\Gamma_1^1$ as in Lemma
\ref{4.K.1} $(i)$.

Regarding the remaining region $\Omega_n/\Gamma$,
we have  $x-y \geq 2-2\delta>0$, $(x,y)\in \Omega_n/\Gamma$.
Moreover, by (\ref{gi-gi}) and the asymptotic of $\gamma_n$ in
\cite[$(2.11)$]{DKMVZ99},
we have $\frac{\gamma_{n-1}^{(n)}}{\gamma_n^{(n)}}=\frac{1}{2}+\mathcal{O}(\frac{1}{n^2})$,
which by the Christoffel-Darboux identity (\cite[3.48]{D99}) implies
\begin{align*} \label{VIn*}
    \mathscr{K}_n^2(x,y)
    \leq & C
           [ (p_n(x;n)p_{n-1}(y;n) )^2
           + ( p_n(y;n)p_{n-1}(x;n))^2 ] e^{-n (V_n(x)+V_n(y))}.
\end{align*}
Hence, in view of the  estimates of orthogonal polynomials  in  Lemma \ref{4.B}
we obtain
\begin{equation} \label{VIn***}
    \iint\limits_{\Omega_n/\Gamma}\mathscr{K}_n^2(x,y)dxdy=\mathcal
    {O}(1).
\end{equation}

Therefore, plugging \eqref{VIn**} and \eqref{VIn***} into
\eqref{VIn*} we prove
\eqref{VIn}.  \hfill $\square$

{\it \bf Proof of Theorem \ref{Bulk-Thm}.}
By virtue of   Propositions \ref{4.exp} and \ref{VIn-est} above,
we can  prove Theorem \ref{Bulk-Thm}
by using similar arguments as in the proof of \cite[Theorems $1.1$ and $1.3$]{G05} or
\cite[Theorems $2.9$ and $2.10$]{Z14.2}.

$(i)$. Take $t$,  $\xi$, $a_n$
and $I_n$ as in Proposition \ref{4.exp}.
Propositions \ref{4.exp} and  \ref{VIn-est} yield
\begin{align} \label{bulk-exp}
   \mathbb{P}_n(\frac{x_k-t}{a_n}<\xi)
   =\mathbb{P}_n(\# I_n \leq n-k)
   =\mathbb{P}_n\(\frac{\# I_n - \mathbb{E}_n \# I_n}{\sqrt{Var_n (\#
   I_n)}}\leq \xi +o(1)\),
\end{align}
which implies the  assertion by the Costin-Lebowitz-Soshnikov theorem (see \cite[p.497-498]{S00}).

$(ii)$.
The proof is based on
the Soshnikov central limit theorem
in \cite[p.174]{S00.2}.
The computations are straightforward but quite complicated.
As in the proof of Proposition
\ref{VIn-est},
the keypoint to calculate the correlation coefficients
$\Lambda_{i,j}$
is that,
similarly to \eqref{VIn**},
for any given
subset $\Lambda \subseteq \Omega_n$ with $\Omega_n$ as in the proof of Proposition \ref{VIn-est},
\begin{align*}
   \iint\limits_{\Lambda} \mathscr{K}^2_n(x,y) dxdy
   =\iint\limits_{\Lambda \cap \wt \Gamma_1^1} \frac{1}{2\pi^2(x-y)^2}
   dxdy + \mathcal{O}(\log\log n),
\end{align*}
where $\wt \Gamma_1^1 = \{(x,y): t \leq x\leq t+\frac{1}{\log n}, t- \frac{1}{\log n} \leq y \leq t-\frac 1n\}$.
For simplicity of exposition,
we refer to \cite{G05} and \cite[Subsection 2.2.2]{Z14.2} for more details.   \hfill $\square$

\subsection{Edge case.} \label{Edge-Section}

We start with Propositions \ref{Edge-Exp} and \ref{Edge-Var} below
concerning the estimates of
expectation and variance, respectively.

\begin{proposition} \label{Edge-Exp}
Set $I:=[t,\9)$ with $t\to 1^-$.
Let  $\wt a(t) := (1-t)^{-\frac 32} \int_t^1 \rho_{V_n}(x)dx$
and $a_1$ be as in Theorem \ref{Edge-Thm}.
Then, we have
\begin{align} \label{asy-edge-exp}
  \mathbb{E}_n(\#I)
    = \wt{a}(t) n(1-t)^{\frac{3}{2}} + \mathcal{O}(1),
\end{align}
and
\begin{align} \label{an-asym}
   \wt a(t) = a_1 + \eta(t) + \calo(n^{-\frac{1}{2m}}),
\end{align}
where
$\eta(t) = -(1-t)^{-3/2} \int_t^1 \wt h'_n(\xi_x)(1-x)^{3/2} dx$
for some $\xi_x\in (x,1)$.
In particular,
for any $t_1, t_2 \in [1-\delta, 1]$ with $\delta \in (0,1)$ fixed,
$|\eta(t_1) - \eta(t_2)| = \calo(|t_1 - t_2|)$.
\end{proposition}

\begin{remark} \label{Rem-Exp-edge}
Let $t$ be the real number such that
$t = 1- (\frac{k}{\wt{a}(t) n})^\frac 23
+ \frac{\sqrt{2}}{3\pi (\wt{a}(t))^\frac 23} \frac{\sqrt{\log k}}{n^\frac 23 k^\frac 13} \xi$,
where $k \sim n^{\g}$, $0<\g<1$.
Then, $\eta(t) = \calo(1-t) = \calo(n^{\frac 23 \g - \frac 23})$.
It follows that for $0<\g\leq \min\{\frac 25, \frac{1}{2m}\}$,
we have
$\mathbb{E}_n(\#I)= a_1  n(1-t)^{\frac{3}{2}} + \mathcal{O}(1)$.
\end{remark}

{\it \bf Proof.}
Similarly to \eqref{compute-Exp}, we have
\begin{align*}
    \mathbb{E}_n(\#I)
    =&\int_t^1 n\rho_{V_n}(x)dx +\mathcal{O}(1).
\end{align*}
Since $t\to 1^-$,
using the definition of $\wt{a}(t)$
we obtain \eqref{asy-edge-exp}.

Moreover,
set $\wt h_n = \frac{1}{2\pi} \sqrt{1+x} h_n(x)$,
where $h_n$ is as in \eqref{rhovn}.
Then,
$\rho_{V_n}(x) = \sqrt{1-x} \wt h_n(x)$.
By Taylor's expansion, mean value theorem and \eqref{4.F.1.1} below,
\begin{align*}
   \int_t^1 \rho_{V_n}(x) dx
   =& \int_t^1 \sqrt{1-x}(\wt h_n(1) + \wt h_n'(\xi_x)(x-1)) dx \\
   =& \int_t^1 \sqrt{1-x} (\frac{\sqrt{2}}{\pi} \sum_{k=0}^{m-1} \frac{A_{m-1-k}}{A_m} + \calo(n^{-\frac{1}{2m}}) - \wt h_n'(\xi_x)(1-x))  dx \\
   =& (1-t)^\frac 32 (a_1  + \calo(n^{-\frac{1}{2m}}) + \eta(t)),
\end{align*}
where $\xi_x\in (x,1)$,
thereby implying \eqref{an-asym}.
The last estimate follows from  the uniform boundedness
of differential $\eta'$ implied by \eqref{4.F.1.2}.
\hfill $\square$

\begin{proposition} \label{Edge-Var}
Let $t$ be such that $t \to 1^-$ and $n(1-t)^{\frac{3}{2}} \to
   \9$, and set $I:=[t,\9)$. Then,
\begin{align} \label{asy-edge-var}
    Var_n(\#I) = \frac{1}{2\pi^2} \log \left[n(1-t)^{\frac{3}{2}}\right](1+o(1)).
\end{align}
\end{proposition}

{\it \bf Proof.} By virtue of
Lemma \ref{4.K.1*}, the estimate
\eqref{asy-edge-var} can be proved by using
similar  arguments as in \cite[Lemma $4$]{S06} and \cite[Proposition $2.12$]{Z14.2}.

In fact, similarly to \eqref{VIn*}, we have
\begin{align} \label{v-k-edge.1}
   Var_n(\# I) =&\iint\limits_{\wt \Gamma}\mathscr{K}_n^2(x,y)dxdy+\iint\limits_{\wt \Omega_n/ \wt \Gamma}\mathscr{K}_n^2(x,y)dxdy,
\end{align}
where $\wt \Omega_n=\{(x,y): t\leq x<\9, -\9<y\leq t\}$, $\wt
\Gamma=\{(x,y):t\leq x\leq t+\frac{1-t}{r_n}, t-\frac{1-t}{r_n}\leq
y\leq t-\e\}$ with $\e = (n\sqrt{1-t})^{-1}$ and
$r^{-1}_n=\max\{\sqrt{1-t}, (\log [n(1-t)^{\frac
32}])^{-1}\}$.

Proceeding as in the proof of \cite[$(3.68)$]{S06}, we
get from \eqref{asy-edge-var.1} that
\begin{align} \label{v-k-edge.2}
    \iint\limits_{\wt \Gamma}\mathscr{K}_n^2(x,y)dxdy
    =\frac{1}{2\pi^2} \log [n(1-t)^{\frac
32}] +\mathcal{O}(\log r_n),
\end{align}
which gives the leading term in
\eqref{asy-edge-var}.

Regarding the remaining integration on $\wt \Omega_n/\wt \Gamma$, using
\eqref{asy-edge-var.1}, Lemma \ref{4.B} and Proposition \ref{Edge-Exp}
and
arguing as in \cite{S06}, we have that
\begin{align} \label{v-k-edge.3}
    \iint\limits_{\wt \Omega_n/\wt \Gamma}\mathscr{K}_n^2(x,y)dxdy
    =\mathcal{O}(\log r_n).
\end{align}
(See also the proof of \cite[(2.3.15)]{Z14.2} for
details.)

Consequently, plugging \eqref{v-k-edge.2} and \eqref{v-k-edge.3}
into \eqref{v-k-edge.1} we obtain  \eqref{asy-edge-var}. \hfill $\square$

{\it \bf Proof of Theorem \ref{Edge-Thm}.}
$(i)$.
Fix $\xi\in \bbr$,
let  $\wt a(t)$ and $t$ be as in Proposition \ref{Edge-Exp} and Remark \ref{Rem-Exp-edge} above,
respectively,
and set $I_n := [t, \9)$, $a_2:= (2\pi^2)^{-\frac 12}$.
 Note that, $t \to 1^-$, as $n \to \9$.
Moreover,
\begin{align*}
    1-t= \(\frac{k}{\wt a(t) n}\)^{\frac 23} \(1- \frac{2a_2 \sqrt{\log k}}{3 k} \xi\),
\end{align*}
which implies that
\begin{align*}
    n(1-t)^{\frac 32}
    = \frac {k}{\wt a(t)} \(1-\frac{2a_2 \sqrt{\log k}}{3 k}\xi \)^{\frac 32}
    =  \frac {k}{\wt a(t)} \(1-\frac{a_2 \sqrt{\log k}}{k} \xi + \mathcal{O}(\frac{\log k}{k^2})\)
    \to \9.
\end{align*}
Then,
define $\wt Y_n$ as in \eqref{def-Yn*}.
Similarly to \eqref{bulk-exp}, by Propositions
\ref{Edge-Exp} and \ref{Edge-Var},
\begin{align} \label{Var-Yn-In}
   \mathbb{P}_n (\wt Y_n <\xi)
   = \mathbb{P}_n(\# I_n \leq k)
   =  \mathbb{P}_n \(\frac{\# I_n - \mathbb{E}_n \# I_n}{\sqrt{Var_n \# I_n}} \leq \xi + o(1)\),
\end{align}
which along with the Costin-Lebowitz-Soshnikov theorem (cf. \cite[p.497-498]{S00}) implies the asymptotic normality of $\wt Y_n$.

Therefore,
taking into account $\wt a(t) \to a_1$ as $n\to \9$,
we see that for any $\tau\in \mathbb{R}$,
$|\mathbb{E}_n e^{i\tau Y_n} - \mathbb{E}_n e^{i\tau \wt Y_n}| \to 0$ as $n\to \9$,
which yields that $Y_n$ is also normally distributed in the limit.

$(ii)$.
The proof is similar but more involved.
Let $\wt s_i$ be such that
$\wt s_i = c_i + d_i \xi_i$,
where
$c_i:=1-(\frac{k_i}{\wt{a}(\wt s_i) n})^{\frac{2}{3}}$
and
$d_i:= \frac{2a_2}{3(\wt{a}(\wt s_i))^\frac 23} \frac{\sqrt{\log k_i}}{n^{\frac{2}{3}}k_i^{\frac{1}{3}}}$
with $a_2 = (2\pi^2 )^{-1/2}$,
$1\leq i\leq m$.
Define $\wt Y_{i,n}$ as in \eqref{def-Yn*}
with $\wt a(\wt s_i)$ replacing $\wt a(t)$,
and set
$Z_i:=\sum_{k=1}^i \# I_k\ (= \#[\wt s_i,\9))$,
$1\leq i\leq m$,
where
$I_1: =[\wt s_1,\9)$
and $I_i=[\wt s_i,\wt s_{i-1})$,
$2\leq i\leq m$.

Since $\wt a(\wt s_i) \to a_1$, $1\leq i\leq m$,
we see that for any $\{c_i\}_{i=1}^m \subseteq \mathbb{R}$
and $\tau \in \mathbb{R}$,
\begin{align*}
  \bigg| \bbe_n e^{i\tau \sum\limits_{i=1}^m c_i Y_{i,n}}
   -  \bbe_n e^{i\tau \sum\limits_{i=1}^m c_i \wt Y_{i,n}} \bigg|
  \to 0, \ \ as\ n\to \9.
\end{align*}
Hence,
it is equivalent to prove the assertion for $\{\wt Y_{i,n}\}_{i=1}^m$.
Moreover,
similarly to \eqref{Var-Yn-In} above, we have
\begin{align*}
   \bbp_n (\wt Y_{1,n}<\xi_1, \cdots, \wt Y_{m,n}<\xi_m)
  = \bbp_n (Z_{1} \leq \xi_1+o(1), \cdots, Z_{m} \leq \xi_m+o(1)).
\end{align*}
Thus,
we can reduce the proof of $\{\wt Y_{i,n}\}_{i=1}^m$
to that of $\{Z_{i}\}_{i=1}^m$.

Now,
the asymptotic normality of $\{Z_{i}\}_{i=1}^m$ can  be proved by using Soshnikov's
central limit theorem (cf. \cite[p.174]{S00.2}),
and the correlations can be calculated by using similar arguments as in the proof
of Proposition \ref{Edge-Var}.
The key fact is that,
for any given set $\Lambda$ in the neighborhood of $(t,t)$ with
$t\to 1^-$ and $n(1-t)^{\frac 32}\to \9$,
\begin{align} \label{edge-corr.2.1*}
    \iint\limits_{\Lambda} \mathscr{K}_n^2(x,y)
    =\frac{1}{2\pi^2} \iint\limits_{\Lambda\cap\ \wt \Gamma}
    \frac{1}{(x-y)^2} dxdy + \mathcal{O}(\log r_n),
\end{align}
where $\wt \Gamma$ and $r_n$ are as in the proof of Proposition
\ref{Edge-Var}.
For more details we refer to the proof of \cite[Theorem 2.14]{Z14.2}
where the constant $a_1$ needs to be modified by the function $\wt a(\cdot)$ here,
yet the arguments still go through since $\wt a(t) = a_1 (1+o(1))$.

It should be mentioned that,
the additional condition $\theta >  \g- \frac{1}{2m}$
arises in the estimate below
(see, e.g., \cite[p.61]{Z14.2})
\begin{align} \label{Dsij-edge}
   \Delta \wt s_{i,j}
   := \wt s_i - \wt s_j
   =\frac{2}{3a_1^{\frac 23}} \frac{1}{n^{\frac 23 + \frac 13 \gamma - \theta}}(1+o(1)),
\end{align}
for which we have to take into account
the delicate remaining order of $\wt a$ in \eqref{an-asym}.
Precisely,
we write $\wt{a}(\wt s_i) = a_1(1+\ve_i)$.
Then,  by the definitions of $\wt s_i$ and $\wt s_j$,
\begin{align*}
  \Delta \wt s_{i,j}
  = \frac{1}{a_1^\frac 23 n^\frac 23}
     \frac{k_j^\frac 23(1+\ve_i)^\frac 23 - k_i^\frac 23(1+\ve_j)^\frac 23}
           {(1+\ve_i)^\frac 23(1+\ve_j)^\frac 23}
     + \calo(\frac{\sqrt{\log k} }{n^{\frac 23 + \frac 13 \g}}).
\end{align*}
Note that $k_j^{\frac 23} - k_i^{\frac 23} =
\frac 23 n^{-\frac 13 \gamma + \theta }(1+o(1))$,
and by \eqref{an-asym},
$\ve_i - \ve_j = \calo(\Delta \wt s_{i,j} + n^{-\frac{1}{2m}})$.
Hence, we get
\begin{align*}
    \Delta \wt s_{i,j}
    = \frac{1}{a_1^\frac 23 n^\frac 23}
      (\frac 23 n^{-\frac 13 \g +\theta}
                     + \calo(n^{\frac 23 \g - \frac {1}{2m}}))  (1+o(1))
                     + o(\Delta \wt s_{i,j}),
\end{align*}
which yields \eqref{Dsij-edge}
under the condition $\theta > \g- \frac{1}{2m}$.
\hfill
$\square$

\begin{remark}
The reason we first treat $\wt Y_n$
(instead of $Y_n$)
in the proof of Theorem \ref{Edge-Thm} is, that
the identity \eqref{Var-Yn-In} may be not valid for $Y_n$
due to the remaining terms of $\wt a(t)$.
Precisely,
if we consider $Y_n$ instead, then
$\mathbb{P}_n (Y_n <\xi)
   =  \mathbb{P}_n \(\frac{\# I_n - \mathbb{E}_n \# I_n}{\sqrt{Var_n \# I_n}} \leq
       \frac{k - \mathbb{E}_n \# I_n}{\sqrt{Var_n \# I_n}}\)$,
where $I_n = [t,\9)$ with
a different choice of $t$ such that
$t = 1- (\frac{k}{a_1 n})^\frac 23
+ \frac{\sqrt{2}}{3\pi a_1^\frac 23} \frac{\sqrt{\log k}}{n^\frac 23 k^\frac 13} \xi$ instead.
However,  Proposition \ref{Edge-Exp} yields
\begin{align*}
   k - \mathbb{E}_n \# I_n
   =& k - \frac{\wt a(t)}{a_1} k  \(1-\frac{a_2 \sqrt{\log k}}{k} \xi + \mathcal{O}(\frac{\log k}{k^2})\) \\
   =&a_2 \sqrt{\log k} (1+o(1))
     + \calo(k\eta(t)) + \calo(kn^{-\frac{1}{2m}}) \\
   =& a_2 \sqrt{\log k} (1+o(1))
      + \calo(n^{\frac 53 \g - \frac 23}) + \calo(n^{\g-\frac{1}{2m}}),
\end{align*}
which implies that
the right-hand side above may exceed the order of $\sqrt{Var_n(\# I_n)}$
(i.e.  $\sqrt{\log k}$)
if $\g> \frac 25$ or $\g> \frac{1}{2m}$.
\end{remark}

\section{Moderate deviations} \label{MDP}

This section contains the proof of Theorem \ref{MDP-Thm} for
unitary invariant ensembles with   Freud-type potentials as in \eqref{v.2}.

{\it \bf Proof of Theorem \ref{MDP-Thm}. }
$(i)$ Fix $\xi \in \mathbb{R}$. Let $t (=t(k,n))$, $a_n$, $t_n$ and $I_n$ be as in  Proposition \ref{4.exp}.
Since $1\ll \g_n \ll \sqrt{\log n}$ and $|t(k,n)|<1$ for $n$ large enough,
we have that $a_n \g_n \xi = o(1)$  and
$|t_n|<1$ for large $n$. Thus, arguing as in \eqref{exp-kn}-\eqref{compute-Exp}, we have
\begin{align} \label{MDP-bulk-exp}
    \mathbb{E}_n(\# I_n)=&n-k- \frac{\sqrt{\log n}}{\sqrt{2\pi^2}}  \g_n \xi+\mathcal
    {O}(1).
\end{align}
Moreover, by Proposition \ref{VIn-est},
\begin{align}   \label{MDP-var-exp}
      Var_n(\# I_n)= & \frac{1}{2\pi^2}\log n
+\mathcal{O}(\log\log n).
\end{align}

Hence, we have
\begin{align}
    \frac{n-k-\mathbb{E}_n(\#I_n)}{ \g_n \sqrt{Var_n(\# I_n)}}
     = \xi + o(1).
\end{align}
This yields that
\begin{align}
   \mathbb{P}_n(\g_n^{-1} X_n <\xi)
   =\mathbb{P}_n(\# I_n \leq n-k)
   =\mathbb{P}_n\(\frac{\# I_n - \mathbb{E}_n \# I_n}{\g_n \sqrt{Var_n(\#
   I_n)}}\leq \xi +o(1)\),
\end{align}
which implies by \cite[Theorem 1.4]{DE13} that for
every $\xi<0$,
\begin{align} \label{bulk-log.1}
     \lim\limits_{n\to \9} \g_n^{-2} \log  \mathbb{P}_n(\g_n^{-1} X_n \leq \xi) = - \frac{\xi^2}{2}.
\end{align}

Similarly, for every $\xi >0$, by \eqref{MDP-bulk-exp} and \eqref{MDP-var-exp},
\begin{align*}
     \mathbb{P}_n(\g_n^{-1} X_n \geq \xi) = \mathbb{P}_n(\# I_n \geq n-k+1)
     =  \mathbb{P}_n\(\frac{\# I_n - \mathbb{E}_n \# I_n}{\g_n \sqrt{Var_n(\#
   I_n)}} \geq \xi + o(1)\),
\end{align*}
which implies by \cite[Theorem 1.4]{DE13} that
\begin{align} \label{bulk-log.2}
   \lim\limits_{n\to \9} \g_n^{-2} \log \mathbb{P}_n (\g_n^{-1} X_n \geq \xi) = - \frac{\xi^2}{2}.
\end{align}

Now, as in the proof of \cite[Theorem 2.1]{DE13}, we denote by $\mathcal{U}$ the set of all open intervals $(c,d)$,
where $c,d\not =0$ and at least one of the endpoints is finite.
Define $\mathcal{L}_{U}:= -\lim_{n\to \9} \g_n^{-2} \log \mathbb{P}_n(\g_n^{-1} X_n \in U)$, $U\in \mathcal{U}$.
By \eqref{bulk-log.1} and \eqref{bulk-log.2},
\begin{align} \label{Lu}
    \mathcal{L}_{U} = \left\{
                        \begin{array}{ll}
                          d^2/2, & \hbox{$c<d<0$;} \\
                          0, & \hbox{$c<0<d$;} \\
                          c^2/2, & \hbox{$0<c<d$.}
                        \end{array}
                      \right.
\end{align}
Then, it follows from \cite[Theorem 4.1.11]{DZ98} that $\{\g_n^{-1}X_n\}$ satisfies a weak  LDP with speed $\g^2_n$ and rate function
$I(x):= \sup_{U\in \mathcal{U},x\in U} \mathcal{L}_U = x^2/2$.

Moreover, for any $\a<\9$, consider the compact set $K_\a = [-c,c]$ with $c = \sqrt{2\a}$.
By \cite[Lemma 1.2.15]{DZ98}
and \eqref{Lu},
\begin{align*}
    &\lim\limits_{n\to \9} \g_n^{-2} \log \mathbb{P}_n(\g_n^{-1} X_n \not \in K_\a) \\
   =& \max\{ \lim\limits_{n\to \9} \g_n^{-2} \log \mathbb{P}_n(\g_n^{-1} X_n \in (-\9,-c)),
       \lim\limits_{n\to \9} \g_n^{-2} \log \mathbb{P}_n(\g_n^{-1} X_n \in (c,\9))\} \\
   =& - \frac{c^2 }{2} = -\a,
\end{align*}
which implies the exponential tightness of $\{\g_n^{-1}X_n\}$, thereby yielding that $\{\g_n^{-1} X_n\}$ satisfies the $LDP(\g_n^2, x^2/2)$.

$(ii)$.
Let $t$ be the real number such that
$t= 1-(\frac{k}{\wt{a}(t)n})^{\frac{2}{3}}
+ \frac{2a_2}{3(\wt{a}(t))^{\frac 23}} \frac{ \sqrt{\log k}}{n^{\frac{2}{3}}k^{\frac{1}{3}}} \g_n  \xi $,
$I_n := [t_n,\9)$,
where $ k, \wt{a}(t)$ are as in Theorem \ref{MDP-Thm},
and $a_2 := (2\pi^2)^{-\frac 12}$.
Since $k/n \to 0$ and $\g_n \ll \sqrt{\log k}$, we have that $t \to 1^{-}$ and
\begin{align*}
  n(1-t)^{\frac 32}
  = \frac{k}{\wt{a}(t)} \(1-\frac{a_2 \sqrt{\log k}}{k} \g_n \xi + \mathcal{O}(\frac{\log k \g_n^2}{k^2})\) \to \9.
\end{align*}

Then,
using Proposition \ref{Edge-Exp} we have
\begin{align}
   \mathbb{E}_n(\#I_n)
   = k-a_2 \sqrt{\log k} \g_n \xi + \calo(\frac{\log k \g_n^2}{k}),
\end{align}
which implies that
\begin{align}
   k - \bbe_n(\#I_n)
   = a_2 \sqrt{\log k} \g_n \xi (1+o(1)).
\end{align}

Moreover, by Proposition  \ref{Edge-Var},
\begin{align}
    \sqrt{Var_n(\# I_n)} = a_2 \sqrt{\log k (1+o(1))}.
\end{align}

Thus, we obtain from the estimates above that
\begin{align*}
    \frac{k - \bbe_n \#I_n}{\g_n \sqrt{Var_n(\# I_n)}} = \xi + o(1),
\end{align*}
which yields that
\begin{align*}
   \mathbb{P}_n(\g_n^{-1} Y_n < \xi)
    =\mathbb{P}_n(\# I_n \leq k)
    = \mathbb{P}_n \( \frac{\# I_n - \mathbb{E}_n \# I_n}{\g_n \sqrt{Var_n(\#
   I_n)}} \leq \xi + o(1)\),
\end{align*}
and
\begin{align*}
   \mathbb{P}_n(\g_n^{-1} Y_n \geq  \xi)
   =\mathbb{P}_n(\# I_n \geq k+1)
    = \mathbb{P}_n \( \frac{\# I_n - \mathbb{E}_n \# I_n}{\g_n \sqrt{Var_n(\#
   I_n)}} \geq \xi + o(1)\).
\end{align*}
Using \cite[Theorem 1.4]{DE13} again we obtain \eqref{bulk-log.1} and \eqref{bulk-log.2}
with $Y_n$ replacing $X_n$.

Therefore,
using similar arguments as those below \eqref{bulk-log.2}
we obtain that $\{\g_n^{-1} Y_n\}$ satisfies the $LDP(\g_n^2, x^2/2)$.
The proof of Theorem \ref{MDP-Thm} is complete. \hfill $\square$

\section{Uniform convex potential} \label{Sec-UC}
This section is devoted to  unitary invariant ensembles with  uniform convex potentials.
Since the arguments are similar to those in the case of Freud-type potentials,
we mainly show the estimates of Christoffel-Darboux kernels
and  orthogonal
polynomials.
Some technical details are
postponed to  the Appendix.

\begin{lemma} \label{rhov}
Consider the uniform convex potential as in \eqref{v.3}. We have
\begin{equation}\label{rhov.1}
    \rho_V(x)=\frac{1}{2\pi}\sqrt{1 -x^2}h(x) \chi_{[-1,1]}(x),
\end{equation}
where $h(x)$ is an analytic function
satisfying
for some $c>0$,
$h(x)\geq c$,  $ \forall x\in \mathbb{R}$.
\end{lemma}
(See the Appendix for the proof.)
As a consequence,  we
have
\begin{lemma} \label{F3}
Consider the uniform convex potential as in \eqref{v.3}.
Let $\delta\in(0,1)$.
Then, $\rho_{V}^{-1}$
and $|\rho'_{V}|$ are  bounded on
$[-1+\delta,1-\delta]$.
\end{lemma}

Set
\begin{align} \label{3.F1}
   F(x):=\left|\int_x^1\frac{1}{2\pi} \sqrt{|1-y^2|}h(y)dy\right|, \ \
   \widetilde{F}(x) :=\left|\int_{-1}^x\frac{1}{2\pi}
   \sqrt{|1-y^2|}h(y)dy\right|.
\end{align}

We also keep the notations  $p_j(x;n)$, $\gamma_j^{(n)}$ and
$\mathscr{K}_n(x,y)$ for the  orthogonal
polynomials, leading coefficients and reproducing kernels with respect
to  $e^{-nV(x)} $.

Below we recall  the Riemann-Hilbert approach developed in \cite{DKMVZ99.0}.
Let $Y(z)$ be an analytic $2\times2$ matrix valued function, solving the Riemann-Hilbert problem
\begin{eqnarray*}
    Y_+(z)=Y_-(z)\left(\begin{array}{ccc}
1&e^{-nV(z)}\\{ }&1\end{array}\right),\ \ for~z\in\mathbb{R}
\end{eqnarray*}
\begin{equation*}
    Y(z)\left(\begin{array}{cc} z^{-n}&{ }\\{ }&z^n\end{array}\right)=I+\mathcal{O}(\frac{1}{|z|}),\ \ as~|z|\rightarrow\infty.
\end{equation*}
Similarly to (\ref{up}),
we have
\begin{equation} \label{yp}
    Y_{11}(z)=\frac{1}{\gamma_n^{(n)}}p_n(z;n),\ \ Y_{21}(z)=-2\pi
i \gamma_{n-1}^{(n)}p_{n-1}(z;n),
\end{equation}

Set
\begin{equation} \label{g}
    g(z) : =\int\log
(z-s)\psi(s)ds, \ \ z\in\mathbb{C}/(-\infty,1),
\end{equation}
where $\psi(z) = \frac{1}{2\pi i} R^{\frac 12}(z) h(z)$, $z\in \mathbb{C}/[-1,1]$,
$R^{\frac{1}{2}}(z)=(z+1)^{\frac{1}{2}}(z-1)^{\frac{1}{2}}$,
which is analytic in $\mathbb{C}/[-1,1]$ and satisfies
$\sqrt{R(z)}\sim z$ as $z\rightarrow\infty$.
Set
\begin{equation} \label{xi}
    G(z): =-\int_1^z
    (s-1)^{\frac{1}{2}}(s+1)^{\frac{1}{2}}h(s)ds,~z\in \mathbb{C}.
\end{equation}
We have
\begin{equation} \label{gvlxi}
    g=\frac{1}{2}(V+l+G),
\end{equation}
where $l$ is as in \eqref{rho-vl1}.
Let
\begin{equation}\label{m}
    M :=e^{-\frac{n l}{2}\sigma_3}Ye^{-n(g-\frac{l}{2})\sigma_3},
\end{equation}
and
\begin{equation} \label{m1}
    M^{(1)}(z) :=\left\{
        \begin{array}{ll}
          M(z), & \hbox{ outside the lens-shaped region;} \\
          M(z)\left(\begin{array}{cc} 1&0\\-e^{-nG}&1\end{array}\right), & \hbox{in the upper lens region;} \\
          M(z)\left(\begin{array}{cc} 1&0\\e^{nG}&1\end{array}\right), & \hbox{in the lower lens region,}
        \end{array}
      \right.
\end{equation}
where the lens  are as in Section \ref{RHM}.

Below are the paramatrices $M_p$ in  small neighborhoods $U_{\pm
1}$ of endpoints $\pm1$.

$(i)$. For $z\in U_1$,
we set (cf. \cite[(4.76)]{DKMVZ99.0})
\begin{equation} \label{mp1}
    M_p:=B(z)P(\Phi_1(z)),
\end{equation}
Here, $B(z)=\frac{1}{\sqrt{2i}}N\left(
                                          \begin{array}{cc}
                                            i & -i \\
                                            1 & 1 \\
                                          \end{array}
                                        \right)
  (\Phi_1)^{\frac{\sigma_3}{4}}$ with $N$ as in  \eqref{N},
\begin{equation} \label{P}
P(z)=\left\{
  \begin{array}{ll}
    \sqrt{2\pi}e^{-\frac{\pi i}{12}} AI(z)e^{(\frac{2}{3}z^{\frac{3}{2}}-\frac{\pi i}{6})\sigma_3}, & \hbox{for $z\in I$;} \\
    \sqrt{2\pi}e^{-\frac{\pi i}{12}} AI(z)e^{(\frac{2}{3}z^{\frac{3}{2}}-\frac{\pi i}{6})\sigma_3}\left(
                                     \begin{array}{cc}
                                       1 & 0 \\
                                       -e^{\frac{4}{3}z^{\frac{3}{2}}} & 1 \\
                                     \end{array}
                                   \right)
, & \hbox{for $z\in II$;} \\
    \sqrt{2\pi}e^{-\frac{\pi i}{12}} \widetilde{AI}(z)e^{(\frac{2}{3}z^{\frac{3}{2}}-\frac{\pi i}{6})\sigma_3}\left(
                                     \begin{array}{cc}
                                       1 & 0 \\
                                       e^{\frac{4}{3}z^{\frac{3}{2}}} & 1 \\
                                     \end{array}
                                   \right), & \hbox{for $z\in III$;} \\
    \sqrt{2\pi}e^{-\frac{\pi i}{12}} \widetilde{AI}(z)e^{(\frac{2}{3}z^{\frac{3}{2}}-\frac{\pi i}{6})\sigma_3}, & \hbox{for $z\in IV$, }
  \end{array}
\right.
\end{equation}
the regions $I$-$IV$ and the matrices $AI$, $\wt{AI}$ are as in Section \ref{RHM},
and
$\Phi_1(z)=(\frac{3n}{4})^{\frac{2}{3}}(-G)^{\frac{2}{3}}=(\frac{3n}{4})^{\frac{2}{3}}(\int_1^zR^{\frac{1}{2}}(s)h(s)ds)^{\frac{2}{3}}$.

$(ii)$. For $z\in U_{-1}$,
we set (cf. \cite[(4.92)]{DKMVZ99.0})
\begin{equation} \label{mp-1}
    M_p:=\widetilde{B}(z)P(\Phi_{-1}(z))\sigma_3,
\end{equation}
where $\widetilde{B}(z)=\frac{1}{\sqrt{2 i}}N\sigma_3\left(
                                            \begin{array}{cc}
                                              i & -i \\
                                              1 & 1 \\
                                            \end{array}
                                          \right)
  (\Phi_{-1})^{\frac{\sigma_3}{4}}
$,
$\Phi_{-1}(z)=(\frac{3n}{4})^{\frac{2}{3}}(-\int_z^{-1}R^{\frac{1}{2}}(s)h(s)ds)^{\frac{2}{3}}$.\\

Finally, set (cf. \cite[(p. 1377)]{DKMVZ99.0})
\begin{equation}\label{Riii}
    R :=\left\{
        \begin{array}{ll}
           M^{(1)}M_p^{-1}, & \hbox{for $z\in U_1\cup U_{-1}$;} \\
          M^{(1)}N^{-1}, & \hbox{otherwise. }
        \end{array}
      \right.
\end{equation}
We have
similar asymptotic expansions of $R$ as in (\ref{4.R}).

Similarly to Lemma \ref{4.K.0}, we have the crucial asymptotic estimates of $\mathscr{K}_n(x,x)$.

\begin{lemma} \label{5.K.0}
Take any sufficiently small $\delta>0$, we have

$(i)$. For $x\in (-1+\delta,1-\delta)$,
\begin{equation} \label{5.K.1*}
    \mathscr{K}_n(x,x)=n\rho_{V}(x)+\mathcal {O}(1).
\end{equation}

$(ii)$. For $x\in(1-\delta,1+\delta)$,
\begin{align} \label{5.K.2}
    \mathscr{K}_n(x,x)=&\left[\frac{1}{4}\frac{\Phi_1'(x)}{\Phi_1(x)}-\frac{a'(x)}{a(x)}\right]2Ai(\Phi_1(x))Ai'(\Phi_1(x))
    \nonumber \\
    &+\Phi_1'(x)\left[(Ai')^2(\Phi_1(x))-\Phi_1(x)Ai^2(\Phi_1(x))\right]+\mathcal {O}(n^{-\frac{5}{6}}).
\end{align}

$(iii)$. For  $x\in(-1-\delta,-1+\delta)$,
\begin{align} \label{5.K.3}
    \mathscr{K}_n(x,x)=&-\left[\frac{1}{4}\frac{\Phi'_{-1}(x)}{\Phi_{-1}(x)}+\frac{a'(x)}{a(x)}\right]2Ai(\Phi_{-1}(x))Ai'(\Phi_{-1}(x)) \nonumber \\
            &-\Phi'_{-1}(x)\left[(Ai')^2(\Phi_{-1}(x))-\Phi_{-1}(x)Ai^2(\Phi_{-1}(x))\right]+\mathcal
            {O}(n^{-\frac{5}{6}}).
\end{align}

$(iv)$. For $x\in \mathbb{R}/(-1-\delta,1+\delta)$,
\begin{equation} \label{5.K.4}
    \mathscr{K}_n(x,x)
    =\frac{1}{4\pi}\frac{1}{(x-1)(x+1)}e^{nG(x)}+\mathcal {O}(n^{-1}).
\end{equation}
\end{lemma}

{\it Proof.}
First note that, similarly to (\ref{k-u}),
\begin{equation} \label{5.K.5}
    2\pi
    i(x-y)\mathscr{K}_n(x,y)=(1,0)Y(x)^TY(y)^{-T}(0,1)^Te^{-n\frac{V(x)+V(y)}{2}}.
\end{equation}

$(i)$. The estimate \eqref{5.K.1*} follows from the same calculations as in Lemma
\ref{4.K.0} $(i)$, with  $U$,
$T$, $S$, $\xi_n$ replaced by $Y$, $M$, $M^{(1)}$ and $G$, respectively.

$(ii)$. First, for $x\in (1-\delta, 1)$, we note  that
$$ N\left(
                                                            \begin{array}{cc}
                                                              i & -i \\
                                                              1 & 1 \\
                                                            \end{array}
                                                          \right)
   =i \left(
                                                            \begin{array}{cc}
                                                              1 & -1 \\
                                                              -i & -i \\
                                                            \end{array}
                                                          \right)
     \left(
       \begin{array}{cc}
         a^{-1} &  \\
          & a \\
       \end{array}
     \right),$$
and $\frac{2}{3}\Phi_1^{\frac{3}{2}}=-\frac{nG}{2}$.
Then, it follows
from (\ref{mp1}) and (\ref{P}) with $z\in II$ that
\begin{align} \label{mp3.1}
    M_p
    =&\widehat{E}_n \left[AI(\Phi_1)\right]e^{-\frac{\pi
i}{6}\sigma_3} \left(
                                           \begin{array}{cc}
                                             1 & 0 \\
                                             -1 & 1 \\
                                           \end{array}
                                         \right)
e^{-\frac{nG}{2}\sigma_3},
\end{align}
where $\widehat{E}_n= \sqrt{\pi} e^{\frac{\pi i}{6}} \left(
                                                       \begin{array}{cc}
                                                         1 & -1 \\
                                                         -i & -i \\
                                                       \end{array}
                                                     \right)
\left(
  \begin{array}{cc}
    a^{-1}(\Phi_1)^{\frac{1}{4}} & 0 \\
    0 & a(\Phi_1)^{-\frac{1}{4}} \\
  \end{array}
\right). $
By (\ref{Riii}),
\begin{equation} \label{m13.1}
    M^{(1)}=R \widehat{E}_n [AI(\Phi_1)] e^{-\frac{\pi i}{6}\sigma_3} \left(
                                           \begin{array}{cc}
                                             1 & 0 \\
                                             -1 & 1 \\
                                           \end{array}
                                         \right)
e^{-\frac{nG}{2}\sigma_3}.
\end{equation}

Similarly, for $x\in
(1,1+\delta)$, using (\ref{P}) with $z\in I$ we have
\begin{equation}\label{mp3.2}
    M_p=\widehat{E}_n\left[AI(\Phi_1)\right]e^{-\frac{\pi
  i}{6}\sigma_3}e^{-\frac{nG}{2}\sigma_3},
\end{equation}
which implies  that
\begin{equation}\label{m13.2}
    M^{(1)}=R \widehat{E}_n\left[AI(\Phi_1)\right]e^{-\frac{\pi
  i}{6}\sigma_3}e^{-\frac{nG}{2}\sigma_3}.
\end{equation}

Then,
comparing \eqref{m13.1} and \eqref{m13.2}
with (\ref{vns1}) and
(\ref{vns2}) below, we see that the
expressions of $M^{(1)}$ and $S$ are similar
on $(1-\delta,1+\delta)$,
with the only difference that
$\wh{E}_n$, $\Phi_1$, $G$ above
are replaced by $E_n$, $f_n$ and $-2\vf_n$,
respectively.
Thus, using the same arguments as in the proof of Lemma
\ref{4.K.0} (ii) we obtain (\ref{5.K.2}).

$(iii)$. First, for $x\in (-1,-1+\delta)$, note that,
\begin{equation*}
    N\sigma_3 \left(
                \begin{array}{cc}
                  i & -i \\
                  1 & 1 \\
                \end{array}
              \right)
    =i \left(
                \begin{array}{cc}
                  1 & 1 \\
                  i & -i \\
                \end{array}
              \right)
    \left(
                \begin{array}{cc}
                  a & 0 \\
                  0 & a^{-1} \\
                \end{array}
              \right)
    \sigma_3,
\end{equation*}
and $\frac{2}{3}(\Phi_{-1})^{\frac{3}{2}}=n(-\frac{G}{2}+\pi i)$.
Using \eqref{mp-1} and (\ref{P}) with $z\in III$  we have
\begin{equation}\label{mp3.3}
    M_p
    =\widehat{\widetilde{E}}_n \sigma_3 \left[\widetilde{AI}(\Phi_{-1})\right]
 e^{-\frac{\pi i}{6}\sigma_3} \left(
  \begin{array}{cc}
    1 & 0 \\
    1 & 1 \\
  \end{array}
\right) \sigma_3 e^{n(-\frac{G}{2}+\pi i)\sigma_3},
\end{equation}
where $\widehat{\widetilde{E}}_n= \sqrt{\pi} e^{\frac{\pi i}{6}}
\left(
                                                       \begin{array}{cc}
                                                         1 & 1 \\
                                                         i & -i \\
                                                       \end{array}
                                                     \right)
\left(
  \begin{array}{cc}
    a(\Phi_{-1})^{\frac{1}{4}} & 0 \\
    0 & a^{-1} (\Phi_{-1})^{-\frac{1}{4}} \\
  \end{array}
\right). $
Hence, by (\ref{Riii}),
\begin{equation} \label{m13.3}
    M^{(1)}=R \widehat{\widetilde{E}}_n \sigma_3 \left[\widetilde{AI}(\Phi_{-1})\right] e^{-\frac{\pi i}{6}\sigma_3} \left(
                                           \begin{array}{cc}
                                             1 & 0 \\
                                             1 & 1 \\
                                           \end{array}
                                         \right)
\sigma_3 e^{n(-\frac{G}{2}+\pi i)\sigma_3}.
\end{equation}

Similarly, for
$x\in(-1-\delta,-1)$, by (\ref{P}) with $z\in IV$,
\begin{equation}\label{mp3.4}
  M_p=\widehat{\widetilde{E}}_n \sigma_3 \left[\widetilde{AI}(\Phi_{-1})\right] e^{-\frac{\pi
i}{6}\sigma_3} \sigma_3 e^{n(-\frac{G}{2}+\pi i)\sigma_3},
\end{equation}
which yields that
\begin{equation} \label{m13.4}
    M^{(1)}=R \widehat{\widetilde{E}}_n \sigma_3 \left[\widetilde{AI}(\Phi_{-1})\right] e^{-\frac{\pi
i}{6}\sigma_3} \sigma_3 e^{n(-\frac{G}{2}+\pi i)\sigma_3}.
\end{equation}

Thus,
comparing \eqref{m13.3} and \eqref{m13.4}
with (\ref{vns3}) and
(\ref{vns4}) below,
we see that
$M^{(1)}$ and
$S$ enjoy the same expressions on $(-1-\delta, -1+\delta)$,
with the only difference that
$\wh{\wt{E}}_n$, $\Phi_{-1}$, $-\frac{G}{2} + \pi i$
here are replaced by $\wt{E}_n$, $-\wt{f_n}$ and $\wt{\vf}_n$,
respectively.
Hence,
similar arguments as in the proof of Lemma
\ref{4.K.0} $(iii)$ yield (\ref{5.K.3}).

 $(iv)$. The proof is
the same as that in Lemma \ref{4.K.0} $(iv)$. $\hfill$ $\square$

Arguing as in the proof of Lemmas \ref{4.K.1} and \ref{4.K.1*}, we have also the asymptotics of $\mathscr{K}_n(x,y)$
as in \eqref{4k6}, \eqref{k7}, \eqref{k7*} and \eqref{asy-edge-var.1}
with $F$ and $\Phi_1$ replacing $F_n$ and $f_n$, respectively.

Below are the asymptotics of orthogonal polynomials,
which would be also of independent interest.
The proof is postponed to the Appendix.

\begin{lemma} \label{5.K.1}
$(i)$. For $x>1+\delta$,
\begin{align}
   p_n(x;n)e^{-\frac{n}{2}V(x)}=&\frac{1}{\sqrt{4\pi}}e^{-n\pi F(x)}
   \left[\left(\frac{x+1}{x-1}\right)^{\frac{1}{4}}+\left(\frac{x-1}{x+1}\right)^{\frac{1}{4}}+\mathcal
{O}(n^{-1})\right], \label{bpn1}\\
   p_{n-1}(x;n)e^{-\frac{n}{2}V(x)}=&\frac{1}{\sqrt{4\pi}}e^{-n \pi F(x)}
   \left[\left(\frac{x+1}{x-1}\right)^{\frac{1}{4}}-\left(\frac{x-1}{x+1}\right)^{\frac{1}{4}}+\mathcal
{O}(n^{-1})\right], \label{bpn2}
\end{align}
where $F(x)$ is defined as in (\ref{3.F1}).

For $x\in(-\9,-1-\delta)$,
\begin{align}
   p_n(x;n)e^{-\frac{n}{2}V(x)}=&(-1)^n \frac{1}{\sqrt{4\pi}}e^{-n\pi \widetilde{F}(x)}
   \left[\left(\frac{x+1}{x-1}\right)^{\frac{1}{4}}+\left(\frac{x-1}{x+1}\right)^{\frac{1}{4}}+\mathcal
{O}(n^{-1})\right], \label{bpn-1}\\
   p_{n-1}(x;n)e^{-\frac{n}{2}V(x)}=&(-1)^n \frac{1}{\sqrt{4\pi}}e^{-n \pi \widetilde{F}(x)}
   \left[\left(\frac{x+1}{x-1}\right)^{\frac{1}{4}}-\left(\frac{x-1}{x+1}\right)^{\frac{1}{4}}+\mathcal
{O}(n^{-1})\right]. \label{bpn-2}
\end{align}
where $\widetilde{F}(x)$ is defined as in (\ref{3.F1}).

$(ii)$. For $x\in(-1+\delta,1-\delta)$,
\begin{align}
   p_n(x;n)e^{-\frac{n}{2}V(x)}=&\sqrt{\frac{2}{\pi}}\frac{1}{(1-x)^{\frac{1}{4}}(1+x)^{\frac{1}{4}}}
   \left\{\cos\left[\frac{1}{2}\arcsin x-\pi n F(x)\right]+\mathcal{O}(n^{-1})\right\},\label{bpn3}\\
   p_{n-1}(x;n)e^{-\frac{n}{2}V(x)}=&\sqrt{\frac{2}{\pi}}\frac{1}{(1-x)^{\frac{1}{4}}(1+x)^{\frac{1}{4}}}
   \left\{\sin\left[\frac{1}{2}\arcsin x+\pi n
   F(x)\right]+\mathcal{O}(n^{-1})\right\}.\label{bpn4}
\end{align}

$(iii)$. For  $x\in (1-\delta,1+\delta)$,
\begin{align}
   p_n(x;n)e^{-\frac{n}{2}V(x)}=&a^{-1}\Phi_{1}^{\frac{1}{4}}Ai(\Phi_1)(1+\mathcal{O}(n^{-1}))
   -a\Phi_1^{-\frac{1}{4}}Ai'(\Phi_1)(1+\mathcal{O}(n^{-1})), \label{bpn5}\\
   p_{n-1}(x;n)e^{-\frac{n}{2}V(x)}=&a^{-1}\Phi_1^{\frac{1}{4}}Ai(\Phi_1)(1+\mathcal{O}(n^{-1}))
   +a\Phi_1^{-\frac{1}{4}}Ai'(\Phi_1)(1+\mathcal{O}(n^{-1})),\label{bpn6}
\end{align}
and for $x\in (-1-\delta,-1+\delta)$,
\begin{align}
   &p_n(x;n)e^{-\frac{n}{2}V(x)} \nonumber \\
   =&(-1)^n\left\{a\Phi_{-1}^{\frac{1}{4}}Ai(\Phi_{-1})(1+\mathcal{O}(n^{-1}))
   -a^{-1}\Phi_{-1}^{-\frac{1}{4}}Ai'(\Phi_{-1})(1+\mathcal{O}(n^{-1}))\right\}, \label{bpn7}\\
   &p_{n-1}(x;n)e^{-\frac{n}{2}V(x)} \nonumber \\
   =&(-1)^{n+1}\left\{a\Phi_{-1}^{\frac{1}{4}}Ai(\Phi_{-1})(1+\mathcal{O}(n^{-1}))
   +a^{-1}\Phi_{-1}^{-\frac{1}{4}}Ai'(\Phi_{-1})(1+\mathcal{O}(n^{-1}))\right\}. \label{bpn8}
\end{align}
\end{lemma}

As a consequence of Lemma \ref{5.K.1}
and the asymptotics of Airy functions,
we also have similar  estimates of  orthogonal polynomials
as in Lemma \ref{4.B}.

Now, by virtue of   the asymptotics of  Christoffel-Darboux kernels and   orthogonal polynomials,
Theorem \ref{Bulk-Thm}, \ref{Edge-Thm} and \ref{MDP-Thm}
in the case of uniform convex potentials can be proved by using similar arguments as in Sections \ref{Gauss-Fluct} and \ref{MDP}.
Note that,
the condition $\theta> \g-1/(2m)$ is not needed here,
since we do not have the  addition order $\calo(n^{-1/(2m)})$  in \eqref{an-asym}.
For simplicity of exposition,
the details are omitted.

\section{Appendix} \label{App}

{\it \bf Proof of Lemma \ref{4.F}.} First, using Theorem
\ref{4.rhon}, we have  for all $n\geq N$ and
$x\in[-1+\delta,1-\delta]$,
\begin{equation*}
    \rho_{V_n}^{-1}(x)=2\pi \frac{1}{\sqrt{1-x^2}}\frac{1}{h_n(x)}
    <\frac{2\pi}{h_0} \frac{1}{\sqrt{1-(1-\delta)^2}}
    <\infty.
\end{equation*}

As regards $\rho'_V$, we have
\begin{equation}\label{4.F.1}
   2\pi
|\rho'_{V_n}|\leq\frac{|x|}{\sqrt{1-x^2}}|h_n(x)|+\sqrt{1-x^2}|h_n'(x)|.
\end{equation}
Note that, by (\ref{vnk}) and (\ref{hn}),
\begin{align}
    h_n(x) =&\sum_{k=0}^{m-1}2\frac{A_{m-k-1}}{A_m}x^{2k}+\mathcal{O}(n^{-\frac{1}{2m}}), \label{4.F.1.1} \\
    h'_n(x) =&\sum_{k=0}^{m-1}4k\frac{A_{m-k-1}}{A_m}x^{2k-1}+\mathcal{O}(n^{-\frac{1}{2m}}),\label{4.F.1.2}
\end{align}
which implies that  $|h_n(x)|$ and $|h'_n(x)|$ are uniformly bounded
for all $n\geq N$ and $x\in[-1+\delta,1-\delta]$,
thereby completing the proof.  \hfill $\square$

{\it \bf Proof of (\ref{kn-en}).} First consider  $x,y \in
(1-\delta,1)$. As in the case where $x,y\in (-1+\delta, 1-\delta)$,
$\eqref{kn+error1***}$ still holds, i.e.,
\begin{equation} \label{4.F.1.2*}
    2\pi
    i(x-y)\mathscr{K}_n(x,y)
    =(e^{-n\varphi_n(x)},e^{n\varphi_n(x)})S(x)^TS(y)^{-T}(-e^{n\varphi_n(y)},e^{-n\varphi_n(y)})^T.
\end{equation}
Since for
$x\in(1-\delta,1)$, $f_n(x+i\epsilon)$ lies in the region $II$ in
(\ref{Psi}), taking $\epsilon\rightarrow 0$ we obtain
\begin{equation*}
    \Psi^{\sigma}(f_n(x))=[AI(f_n(x))]e^{-\frac{\pi i}{6}\sigma_3}\left(
                                                   \begin{array}{cc}
                                                     1 & 0 \\
                                                     -1 & 1 \\
                                                   \end{array}
                                                 \right),
\end{equation*}
which along with (\ref{R}) and (\ref{Pn1}) yields that
\begin{equation} \label{vns1}
    S=RE_n[AI(f_n)]e^{-\frac{\pi i}{6}\sigma_3}\left(
                                                   \begin{array}{cc}
                                                     1 & 0 \\
                                                     -1 & 1 \\
                                                   \end{array}
                                                 \right)
      e^{n\varphi_n\sigma_3}.
\end{equation}
Consequently, plugging \eqref{vns1} into \eqref{4.F.1.2*},
we obtain (\ref{kn-en}) for
$x,y\in(1-\delta,1)$.

Regarding the case where $x,y\in (1,1+\delta)$, by (\ref{T}), (\ref{S}), (\ref{k-u}) and
(\ref{gn}),
\begin{equation*}
    U=e^{n\frac{l_n}{2}\sigma_3}Se^{n(g_n-\frac{l_n}{2})\sigma_3},
\end{equation*}
and
\begin{align} \label{vns1*}
    2\pi i(x-y)\mathscr{K}_n(x,y)
    =(e^{-n\varphi_n(x)},0)S^T(x)S^{-T}(y)(0,e^{-n\varphi_n(y)})^T.
\end{align}
Since for $x\in(1,1+\delta)$, $f_n(x+i\epsilon)$ is in the region $I$
in (\ref{Psi}), by (\ref{R}) and (\ref{Pn1}),
\begin{equation} \label{vns2}
    S=RE_n[AI(f_n)]e^{-\frac{\pi
    i}{6}\sigma_3}e^{n\varphi_n\sigma_3}.
\end{equation}
Hence, combining \eqref{vns1*} and \eqref{vns2}  we get
(\ref{kn-en}) for
$x,y\in(1,1+\delta)$, thereby completing the proof of \eqref{kn-en}. \hfill $\square$

{\it \bf Proof of (\ref{kn-ai}).} We first show that $I_2(x,y) = \calo(n^{-\frac{5}{6}})$.
Indeed,
by \eqref{I2}, expressions of $AI$,
$E_n$ and that $\det [AI(z)]= \frac{-1}{2\pi i} e^{-\frac{\pi i}{3}}$
(see \cite[p.$890$]{DG07}),
\begin{align*}
    I_2(x,y)
    =&-2\pi i e^{\frac{\pi i}{3}}(H_n(x)Ai(f_n(x)),H_n^{-1}(x)Ai'(f_n(x)))\left(
                                                  \begin{array}{cc}
                                                    1 & -1 \\
                                                    -i & -i \\
                                                  \end{array}
                                                \right)^T
    \\
     & \cdot \Delta_{R}(x,y)R^{-T}(y)\left(
                                                \begin{array}{cc}
                                                    1 & -1 \\
                                                    -i & -i \\
                                                  \end{array}
                                                \right)^{-T}
    (-H_n^{-1}(y)Ai'(f_n(y)),H_n(y)Ai(f_n(y)))^T.
\end{align*}
Note that, since $H_n=f_n^{\frac{1}{4}}a^{-1}$, by (\ref{fn2}),
\begin{equation*}
    H_n
    =n^{\frac{1}{6}}(x-1)^{\frac{1}{4}}(\widehat{\phi}_n)^{\frac{1}{6}}(\frac{x-1}{x+1})^{-\frac{1}{4}}
    =n^{\frac{1}{6}}(x+1)^{\frac{1}{4}}(\widehat{\phi}_n)^{\frac{1}{6}}
    =\mathcal{O}(n^{\frac{1}{6}}).
\end{equation*}
Moreover,  for $x\in\mathbb{R}$,  $|Ai(x)|=\mathcal{O}(1)$ and
$|Ai'(f_n(x))|=\mathcal{O}(|f_n(x)|^{\frac{1}{4}})=\mathcal{O}(n^{\frac{1}{6}})$,
and by (\ref{4.R}),
$\Delta{R}(x,y)R^{-T}(y)=\mathcal{O}(n^{-1})$. Thus,
we conclude that $I_2(x,y)$ is of
order
$\mathcal{O}(n^{\frac{1}{6}})\mathcal{O}(n^{-1})=\mathcal{O}(n^{-\frac{5}{6}})$.

It remains to check the first term on the right-hand side of
(\ref{kn-error2}). To this end, it follows from \eqref{kn-error2} and
the computations as above that
\begin{align*}
    &e^{-\frac{\pi i}{3}}(1,0)[AI(f_n(x))]^TE_n^T(x)E_n^{-T}(y)[AI(f_n(y))]^{-T}(0,1)^T\\
  =&(-2\pi i)\left[-Ai(f_n(x))Ai'(f_n(y))\frac{f_n^{\frac{1}{4}}(x)}{f_n^{\frac{1}{4}}(y)}\frac{a(y)}{a(x)}
   +Ai'(f_n(x))Ai(f_n(y))\frac{f_n^{\frac{1}{4}}(y)}{f_n^{\frac{1}{4}}(x)}\frac{a(x)}{a(y)}\right].
\end{align*}
which yields the first term
on the right-hand side of (\ref{kn-ai}).\hfill $\square$

{\it \bf Proof of (\ref{kn-en2}).} The proofs are similar to those
of (\ref{kn-en}). First consider  $x,y\in(-1,-1+\delta)$. As in
the case where $x,y\in (1-\delta,1)$, we have
\begin{equation} \label{vns3*}
    2\pi i(x-y)\mathscr{K}_n(x,y)
    =(e^{-n\varphi_n(x)},e^{n\varphi_n(x)})S(x)^TS(y)^{-T}(-e^{n\varphi_n(y)},e^{-n\varphi_n(y)})^T.
\end{equation}
Since for $x\in(-1,-1+\delta)$, $-\widetilde{f}_n(x+i\epsilon)$
lies in the region $III$ in (\ref{Psi}), letting
$\epsilon\rightarrow 0$ we have
\begin{equation*}
    \wt \Psi^{\sigma}(-\widetilde{f}_n(x))=[\widetilde{AI}(-\widetilde{f}_n(x))]e^{-\frac{\pi i}{6}\sigma_3}\left(
                                                                                                         \begin{array}{cc}
                                                                                                           1 & 0 \\
                                                                                                           1 & 1 \\
                                                                                                         \end{array}
                                                                                                       \right),
\end{equation*}
which along with (\ref{R}) and (\ref{Pn-1}) yields
\begin{equation} \label{vns3}
    S=R\widetilde{E}_n\sigma_3[\widetilde{AI}(-\widetilde{f}_n)]e^{-\frac{\pi i}{6}\sigma_3}\left(
                                                                                                         \begin{array}{cc}
                                                                                                           1 & 0 \\
                                                                                                           1 & 1 \\
                                                                                                         \end{array}
                                                                                                       \right)
    \sigma_3e^{n\widetilde{\varphi}_n\sigma_3}.
\end{equation}
Thus, plugging \eqref{vns3} into \eqref{vns3*}, since
$\widetilde{\varphi}_n(z)=\varphi_n(z)+\pi i$, $z\in\mathbb{C}^+$,
we get (\ref{kn-en2}).

Regarding the case where $x,y\in(-1-\delta,-1)$. As in the case where
$x,y\in (1,1+\delta)$ in the proof of \eqref{kn-en}, we have
\begin{equation} \label{vns3**}
    2\pi
    i(x-y)\mathscr{K}_n(x,y)=(e^{-n\varphi_n(x)},0)S^T(x)S^{-T}(y)(0,e^{-n\varphi_n(y)})^T.
\end{equation}
Since for $x\in(-1-\delta,-1)$, $-\widetilde{f}_n(x+i\epsilon)$ is
in the region $IV$ in (\ref{Psi}),  taking $\ve \to 0$,  we obtain from \eqref{Pn-1}
and (\ref{R}) that
\begin{equation} \label{vns4}
    S=R\widetilde{E}_n\sigma_3 [\widetilde{AI}(-\widetilde{f}_n)] e^{-\frac{\pi
i}{6}\sigma_3}\sigma_3e^{n\widetilde{\varphi}_n\sigma_3}.
\end{equation}
Therefore, plugging \eqref{vns4} into \eqref{vns3**} yields (\ref{kn-en2}). \hfill $\square$

{\it \bf Proof of Lemma \ref{rhov}}.
Define the
Hilbert transform $\mathscr{H}$ and the Borel transform $\mathscr{B}$ by
\begin{align*}
    \mathscr{H} \rho_V(x)=&\frac{1}{\pi}P.V.\int\frac{\rho_V(y)}{x-y}dy,\\
    \mathscr{B}\rho_V(z)=&\frac{1}{\pi i}\int\frac{\rho_V(s)}{s-z}ds, \  \
    z\in\mathbb{C}/\mathbb{R}.
\end{align*}

In view of \cite[$(3.10)$, $(3.12)$]{DKMVZ99.0}, for  $x\in\mathbb{R}$ we have
\begin{align*}
    (\mathscr{B}\rho_V)_{\pm}(x)=&\pm\rho_V(x)+i\mathscr{H} \rho_V(x)=\pm\rho_V(x)-\frac{1}{2\pi
i}V'(x).
\end{align*}
Moreover, by virtue of \cite[$(3.17)$, $(3.18)$]{DKMVZ99.0},
\begin{equation*}
    \mathscr{B}\rho_V(z)=-\frac{1}{2\pi i}V'(z)-\frac{\sqrt{R(z)}}{4\pi^2
}\oint\limits_{\Gamma_z}\frac{V'(s)}{\sqrt{R(s)}}\frac{ds}{s-z},
\end{equation*}
where $\sqrt{R(z)}$ is as in Section \ref{Sec-UC}, and $\Gamma_z$
is a counterclockwise contour with $z$ and $[-1,1]$ in its interior.
Note that, due to the analytic branch  of
$\sqrt{R(z)}$, we have
\begin{equation} \label{R+-}
    (\sqrt{R(x)})_+=i\sqrt{(x+1)(1-x)}=-(\sqrt{R(x)})_-.
\end{equation}

Hence, it follows that for
$x\in(-1,1)$,
\begin{equation*}
    \rho_V(x)=\frac{1}{2\pi}\sqrt{(1-x)(x+1)} \(\frac{1}{2\pi
i}\oint\limits_{\Gamma_z}\frac{V'(s)}{\sqrt{R(s)}}\frac{ds}{s-z}\)_+ ,
\end{equation*}
and so
\begin{equation}\label{h.1}
    h(z)=\frac{1}{2\pi
i}\oint\limits_{\Gamma_z}\frac{V'(s)}{\sqrt{R(s)}}\frac{ds}{s-z},
\end{equation}
which is an analytic function.

It remains to prove that  for some $c>0$,
$h(x)\geq c$ $\forall x\in\mathbb{R}$.
We first claim that
\begin{equation} \label{h.1*}
    \frac{1}{2\pi
i}\oint\limits_{\Gamma_z}\frac{1}{\sqrt{R(s)}}\frac{ds}{s-z}\equiv0
\end{equation}
To this end, we have for $r$ large enough that
\begin{align} \label{h.1**}
    \frac{1}{2\pi i}\oint\limits_{\Gamma_z}\frac{1}{\sqrt{R(s)}}\frac{ds}{s-z}
    =&\frac{1}{2\pi
i}\oint\limits_{|s|=r}\frac{1}{\sqrt{(s+1)(s-1)}}\frac{ds}{s-z}.
    \end{align}
Using Taylor's extension we see that
\begin{equation*}
    \frac{1}{\sqrt{s^2-1}}\frac{1}{s-z}
    =(A_0\frac{1}{s}+A_1\frac{1}{s^3}+A_2\frac{1}{s^5}+...)
    (\frac{1}{s}+\frac{z}{s^2}+\frac{z^2}{s^3}+...),
\end{equation*}
which by Cauchy's theorem implies that
the right-hand side of \eqref{h.1**} equals to the coefficient of
$\frac{1}{s}$
which is exactly zero, thereby yielding \eqref{h.1*}, as claimed.

Now, it follows from (\ref{h.1}) and \eqref{h.1*} that
\begin{align} \label{h.1***}
    h(z)
    =\frac{1}{2\pi
i}\oint\limits_{\Gamma_z}\frac{V'(s)-V'(z)}{s-z}\frac{1}{\sqrt{R(s)}}ds
    =\frac{1}{2\pi
i}\oint\limits_{\Gamma_1}\frac{V'(s)-V'(z)}{s-z}\frac{1}{\sqrt{R(s)}}ds,
\end{align}
where $\Gamma_1$
is a counterclockwise contour with $[-1,1]$, but not $z$, in the interior.

Therefore, in view of \eqref{h.1***} we have that for $x\in(-1,1)$,
\begin{align*}
   h(x)
   =&\lim\limits_{z\in\mathbb{C}^+\rightarrow x}h(z)
   =\frac{1}{\pi
i}\int_1^{-1}\frac{V'(s)-V'(x)}{s-x}\frac{1}{(\sqrt{R(s)})_+}ds \\
   =& \frac{1}{\pi} \int_{-1}^1\frac{V'(s)-V'(x)}{s-x}\frac{1}{\sqrt{(s+1)(1-s)}}ds,
\end{align*}
which implies by the  mean value theorem and the uniform convexity of $V$ that
\begin{align*}
  h(x)
   =\frac{1}{\pi}\int_{-1}^1\frac{V''(\xi)}{\sqrt{(s+1)(1-s)}}ds
   \geq\frac{c}{\pi}\int_{-1}^1\frac{ds}{\sqrt{(s+1)(1-s)}}
   =&c>0.
\end{align*}
where $\xi\in(-1,1)$. The proof is complete. \hfill $\square$

Before proving Lemma \ref{5.K.1} we recall that
\begin{theorem} (\cite[(1.62) -- (1.64)]{DKMVZ99.0})  \label{gamma3}
Consider the uniform convex potential as in \eqref{v.3}.
For the leading coefficients of  orthogonal polynomials we have
\begin{align*}
   (\gamma_{n-1}^{(n)})^2=e^{-nl}[\frac{1}{4\pi} + \mathcal{O}(n^{-1})], \
   (\gamma_{n}^{(n)})^{-2}=e^{nl} [\pi  + \mathcal{O}(n^{-1})], \
   \frac{\gamma_{n-1}^{(n)}}{\gamma_{n}^{(n)}}=\frac{1}{2}+\mathcal{O}(n^{-1/2}).
\end{align*}
\end{theorem}

\begin{theorem} (\cite[Theorem 1.1 -- 1.3]{DKMVZ99.0}) \label{PR3}
Consider the uniform convex potential as in \eqref{v.3}.
For the monic polynomials  we have

$(i)$. For $x\in \mathbb{R}/(-1-\delta,1+\delta)$,
\begin{align}
\pi_n(x;n)=&e^{ng(x)}\left(M_1(x)+\mathcal{O}(n^{-1})\right),\label{pin-m1}\\
    -2\pi
i(\gamma_{n-1}^{(n)})^2\pi_{n-1}(x;n)=&e^{n(g(x)-l)}\left(M_2(x)+\mathcal{O}(n^{-1})\right),   \label{pin-m2}
\end{align}
where  $M_1=\frac{a+a ^{-1}}{2}$, $M_2=\frac{a^{-1}-a}{2i}$ with  $a$ as in (\ref{a})
and $l$ is as in \eqref{rho-vl1}.

$(ii)$. For $x\in(-1+\delta,1-\delta)$,
\begin{align}
    \pi_n(x;n)=&2e^{\frac{n}{2}(V(x)+l)}\left[Re(M_1e^{i\pi
n F(x)})+\mathcal{O}(n^{-1})\right], \label{pin-sin}\\
    -2\pi i(\gamma_{n-1}^{(n)})^2\pi_{n-1}(x;n)=&2e^{\frac{n}{2}(V(x)-l)}\left[Im(M_2e^{i\pi
n F(x)})+\mathcal{O}(n^{-1})\right],\label{pin-con}
\end{align}
where $F$ is as in (\ref{3.F1}).

$(iii)$. For $x\in (1-\delta,1)\cup (-1,-1+\delta)$,
\begin{align}
    \pi_n(x;n)=&\left(e^{\frac{nl}{2}\sigma_3}(I+\mathcal{O}(n^{-1}))M_pe^{n(g(x)-\frac{l}{2})\sigma_3}\left(
                                                                                                                       \begin{array}{cc}
                                                                                                                         1 & 0 \\
                                                                                                                         e^{nV} & 1 \\
                                                                                                                       \end{array}
                                                                                                                     \right)
\right)_{11},\label{pin-mp3}\\
    -2\pi
i(\gamma_{n-1}^{(n)})^2\pi_{n-1}(x;n)=&\left(e^{\frac{nl}{2}\sigma_3}(1+\mathcal{O}(n^{-1}))M_pe^{n(g(x)-\frac{l}{2})\sigma_3}\left(
                                                                                                                       \begin{array}{cc}
                                                                                                                         1 & 0 \\
                                                                                                                         e^{nV} & 1 \\
                                                                                                                       \end{array}
                                                                                                                     \right)
\right)_{21},\label{pin-mp4}
\end{align}\\
and for $x\in (1,1+\delta)\cup (-1-\delta,-1)$,
\begin{align}
    \pi_n(x;n)=&\left((I+\mathcal{O}(n^{-1}))M_p\right)_{11}e^{ng(x) },\label{pin-mp1}\\
    -2\pi
i(\gamma_{n-1}^{(n)})^2\pi_{n-1}(x;n)=&\left((I+\mathcal{O}(n^{-1}))M_p\right)_{21}e^{n(g(x)-l)},\label{pin-mp2}
\end{align}
\end{theorem}

{\it \bf Proof of Lemma \ref{5.K.1}.} We first note that, by the analytic branch of $R^{1/2}(z)$,
\begin{equation*}
    (z-1)^{\frac{1}{2}}=i(1-z)^{\frac{1}{2}},~(z+1)^{\frac{1}{2}}=i(-1-z)^{\frac{1}{2}}, \ z\in\mathbb{C}_+.
\end{equation*}
Then, using  (\ref{xi}) we have
\begin{align}
   G(x)=&-\int_1^x\sqrt{(s-1)(s+1)}h(s)ds=-2 \pi F(x), \ \ x>1, \label{G-F.1} \\
   G(x)=&2\pi i F(x), \ \ x\in(-1,1),  \label{G-F.2} \\
   G(x)=&- 2\pi \widetilde{F}(x)+2\pi i, \ \ x<-1,  \label{G-F.3}
\end{align}
where $F(x)$ and $\widetilde{F}(x)$ are as in (\ref{3.F1}).

$(i)$. For $x\in \mathbb{R}/(-1-\delta,1+\delta)$, by
(\ref{pin-m1}) and (\ref{gvlxi}),
\begin{equation*}
    \pi_n(x;n)=e^{\frac{n}{2}(V+l+G)}\left\{\frac{1}{2}\left[(\frac{x+1}{x-1})^{\frac{1}{4}}+(\frac{x-1}{x+1})^{\frac{1}{4}}\right]+\mathcal
{O}(n^{-1})\right\}.
\end{equation*}
Then, by Theorem \ref{gamma3},
\begin{equation*}
    p_n(x;n)=\gamma_n^{(n)}\pi_n(x;n)=\frac{1}{\sqrt{\pi}}e^{\frac{n}{2}(V+G)}\left\{\frac{1}{2}\left[(\frac{x+1}{x-1})^{\frac{1}{4}}+(\frac{x-1}{x+1})^{\frac{1}{4}}\right]+\mathcal
{O}(n^{-1})\right\}.
\end{equation*}
Thus, using (\ref{G-F.1}) and (\ref{G-F.3}) for $x > 1+\delta$ and
$x< -1-\delta$, respectively,
we obtain (\ref{bpn1}) and
(\ref{bpn-1}).

Similarly, by Theorem \ref{gamma3}, (\ref{pin-m2}) and
(\ref{gvlxi}),
\begin{align*}
    p_{n-1}(x;n)e^{-\frac{n}{2}V}
     =&\frac{1}{\sqrt{4 \pi}} e^{\frac{ n G}{2}} \left\{\left[(\frac{x+1}{x-1})^{\frac{1}{4}}-(\frac{x-1}{x+1})^{\frac{1}{4}}\right]+\mathcal
{O}(n^{-1})\right\},
\end{align*}
which yields (\ref{bpn2}) and (\ref{bpn-2}) by (\ref{G-F.1}) and
(\ref{G-F.3}), respectively.

$(ii)$. By the definitions of $M_1$ and $M_2$, we have that  (cf.   \cite[(8.33), (8.34)]{DKMVZ99})
\begin{align*}
    M_1=&\frac{a+a^{-1}}{2}=\frac{\sqrt{2}}{2}\frac{1}{(1-x)^{\frac{1}{4}}(1+x)^{\frac{1}{4}}}e^{-\frac{i}{2}\arcsin
x}\\
    M_2=&\frac{a^{-1}-a}{2i}=-\frac{\sqrt{2}}{2}\frac{1}{(1-x)^{\frac{1}{4}}(1+x)^{\frac{1}{4}}}e^{\frac{i}{2}\arcsin
x}.
\end{align*}
Then, in view of Theorem \ref{gamma3} and Theorem \ref{PR3} $(ii)$, we obtain  (\ref{bpn3}) and
(\ref{bpn4}).

$(iii)$. We consider four cases $(iii.1)-(iii.4)$ below.

$(iii.1)$. For $x\in(1-\delta,1)$, by (\ref{mp3.1}) and
(\ref{gvlxi}),
\begin{align*}
     e^{\frac{nl}{2}\sigma_3}M_pe^{n(g-\frac{l}{2})\sigma_3}\left(
                                                                    \begin{array}{cc}
                                                                    1 & 0 \\
                                                                    e^{nV} & 1 \\
                                                                    \end{array}
                                                                    \right)
    =e^{\frac{nl}{2}\sigma_3}\widehat{E}_n \left[AI(\Phi_1)\right]
  e^{-\frac{i\pi}{6}\sigma_3}e^{\frac{n}{2}V\sigma_3}.
\end{align*}
Then, since $\widehat{E}_n=\sqrt{\pi} e^{\frac{\pi i}{6}}\left(
                                                         \begin{array}{cc}
                                                           a^{-1}\Phi_1^{\frac{1}{4}} & -a\Phi_1^{-\frac{1}{4}} \\
                                                           -ia^{-1}\Phi_1^{\frac{1}{4}} & -ia\Phi_1^{-\frac{1}{4}} \\
                                                         \end{array}
                                                       \right)
$, direct calculations show that
\begin{equation} \label{c11}
   \left(e^{\frac{nl}{2}\sigma_3}M_pe^{n(g-\frac{l}{2})\sigma_3}\left(
                                                                    \begin{array}{cc}
                                                                    1 & o \\
                                                                    e^{nV} & 1 \\
                                                                    \end{array}
                                                                    \right)\right)_{11}
=\sqrt{\pi}e^{\frac{n}{2}(V+l)}
   \left[a^{-1}\Phi_1^{\frac{1}{4}}Ai(\Phi_1)-a\Phi_1^{-\frac{1}{4}}Ai'(\Phi_1)\right].
\end{equation}

Similarly,
\begin{equation} \label{c21}
   \left(e^{\frac{nl}{2}\sigma_3}M_pe^{n(g-\frac{l}{2})\sigma_3}\left(
                                                                    \begin{array}{cc}
                                                                    1 & o \\
                                                                    e^{nV} & 1 \\
                                                                    \end{array}
                                                                    \right)\right)_{21}
=(-i) \sqrt{\pi}e^{\frac{n}{2}(V-l)}
   \left[a^{-1}\Phi_1^{\frac{1}{4}}Ai(\Phi_1)+a\Phi_1^{-\frac{1}{4}}Ai'(\Phi_1)\right].
\end{equation}
Plugging these into (\ref{pin-mp3}) and (\ref{pin-mp4})  and
using Theorem \ref{gamma3}, we  obtain
(\ref{bpn5}) and (\ref{bpn6}).

$(iii.2)$. For $x\in (1,1+\delta)$, using (\ref{mp3.2}) we note that in (\ref{pin-mp1}) and
(\ref{pin-mp2}), $(M_p)_{11}e^{ng}$ and $(M_p)_{21}e^{n(g-l)}$ have
the same formulations as (\ref{c11}) and (\ref{c21}). Thus,
arguing as above we obtain (\ref{bpn5}) and (\ref{bpn6}).

$(iii.3)$. For $x\in (-1,-1+\delta)$, it follows from (\ref{mp3.3})
that
\begin{align*}
     e^{\frac{nl}{2}\sigma_3} M_p e^{n(g-\frac{l}{2})\sigma_3}\left(
                                                               \begin{array}{cc}
                                                                 1 & 0 \\
                                                                 e^{nV} & 1 \\
                                                               \end{array}
                                                             \right)
     =(-1)^n e^{\frac{nl}{2}\sigma_3} \widehat{\widetilde{E}}_n \sigma_3 [\widetilde{AI}(\Phi_{-1})] e^{-\frac{\pi
    i}{6}\sigma_3} e^{\frac{n}{2}V\sigma_3} \sigma_3.
\end{align*}
Then, since $\widehat{\widetilde{E}}_n =\sqrt{\pi} e^{\frac{\pi
i}{6}}\left(
                                                                     \begin{array}{cc}
                                                                       a\Phi_{-1}^{\frac{1}{4}} & a^{-1}\Phi_{-1}^{-\frac{1}{4}} \\
                                                                       ia\Phi_{-1}^{\frac{1}{4}} & -ia^{-1}\Phi_{-1}^{-\frac{1}{4}} \\
                                                                     \end{array}
                                                                   \right)
$, we have
\begin{align}\label{d11}
  &\left(e^{\frac{nl}{2}\sigma_3} M_p e^{n(g-\frac{l}{2})\sigma_3}\left(
                                                               \begin{array}{cc}
                                                                 1 & 0 \\
                                                                 e^{nV} & 1 \\
                                                               \end{array}
                                                             \right)\right)_{11}
                                                             \nonumber    \\
  =&(-1)^n \sqrt{\pi} e^{\frac{n(V+l)}{2}} \left[a \Phi_{-1}^{\frac{1}{4}}  Ai(\Phi_{-1}) - a^{-1} \Phi_{-1}^{-\frac{1}{4}}
  Ai'(\Phi_{-1})\right].
\end{align}

Similarly,
\begin{align} \label{d21}
  &\left(e^{\frac{nl}{2}\sigma_3} M_p e^{n(g-\frac{l}{2})\sigma_3}\left(
                                                               \begin{array}{cc}
                                                                 1 & 0 \\
                                                                 e^{nV} & 1 \\
                                                               \end{array}
                                                             \right)\right)_{21}
                                                             \nonumber                     \\
  =&(-1)^{n+1} (-i) \sqrt{\pi} e^{\frac{n(V-l)}{2}} \left[a \Phi_{-1}^{\frac{1}{4}}  Ai(\Phi_{-1}) + a^{-1} \Phi_{-1}^{-\frac{1}{4}}  Ai'(\Phi_{-1})\right].
\end{align}
Thus, (\ref{bpn7}) and (\ref{bpn8}) follow from (\ref{pin-mp3}),
(\ref{pin-mp4}) and Theorem \ref{gamma3}.

$(iii.4)$ For $x\in (-1-\delta,-1)$, by (\ref{mp3.4}) we note that
$(M_p)_{11}e^{ng}$ and $(M_p)_{11}e^{n(g-l)}$ in (\ref{pin-mp1}) and
(\ref{pin-mp2}) have the same formulations as (\ref{d11}) and (\ref{d21}),
which consequently implies  (\ref{bpn7}) and (\ref{bpn8}).
The proof
of Lemma \ref{5.K.1} is complete. \hfill $\square$

{\it \bf Acknowledgement.} The author is grateful to
Prof. Xiang-Dong Li for  valuable discussions
and Prof. Zhonggen Su for  nice lectures at AMSS in Beijing in 2011
when this work was initiated.
The author would also like to thank the two referees
for  valuable comments to improve this paper.
In particular,
I would like to thank one referee for pointing out the reference \cite{BD14}.
Financial support by the NSFC (No. 11871337) is also gratefully
acknowledged.


\begin{thebibliography}{nn}

\bibitem{BEY14} P. Bourgade,  L. Erd\"os, H. T. Yau,
Edge universality of beta ensembles.
Comm. Math. Phys. 332 (2014), no. 1, 261-353.

\bibitem{BD14}
J. Breuer, M. Duits,
The Nevai condition and a local law of large numbers for orthogonal polynomial ensembles.
Adv. Math. 265 (2014), 441-484.

\bibitem{CL95}  O. Costin,  J. L. Lebowitz,
Gaussian fluctuation in random matrices.
Phys. Rev. Lett. 75 (1995), no. 1, 69-72.


\bibitem{C92}  T. Chan,
The Wigner semi-circle law and eigenvalues of matrix-valued diffusions.
Probab. Theory Related Fields 93 (1992), no. 2, 249-272.

\bibitem{D99} P. Deift, Orthogonal polynomials and random matrices: A Riemann-Hilbert approach,
Courant Lecture Notes in Mathematics 3, Courant Institute of
Mathematical Science, New York, 1999.

\bibitem{DG07}  P. Deift, D. Gioev,
Universality at the edge of the spectrum for unitary, orthogonal, and symplectic ensembles of random matrices.
Comm. Pure Appl. Math. 60 (2007), no. 6, 867-910.


\bibitem{DKMVZ99.0}   P. Deift, T. Kriecherbauer, K. T.-R McLaughlin,  S. Venakides, X. Zhou,
Uniform asymptotics for polynomials orthogonal with respect to varying exponential weights and applications
to universality questions in random matrix theory. Comm. Pure Appl. Math. 52 (1999), no. 11, 1335-1425.


\bibitem{DKMVZ99}  P. Deift, T. Kriecherbauer, K. T.-R McLaughlin, S. Venakides, X. Zhou,
Strong asymptotics of orthogonal polynomials with respect to exponential weights.
Comm. Pure Appl. Math. 52 (1999), no. 12, 1491-1552.

\bibitem{DKMVZ01} P. Deift, T. Kriecherbauer, K. T.-R McLaughlin, S. Venakides,
Zhou, X.: A Riemann-Hilbert approach to asymptotic questions for
orthogonal polynomials. J. Comput. Appl. Math. 133 (2001), no. 1-2, 47-63.

\bibitem{DVZ97}  P. Deift, S.  Venakides, X. Zhou,
New results in small dispersion KdV by an extension of the steepest descent method for Riemann-Hilbert problems.
Internat. Math. Res. Notices 1997, no. 6, 286-299.

\bibitem{DZ93}  P. Deift, X. Zhou, A steepest descent method for oscillatory Riemann-Hilbert problems.
Asymptotics for the MKdV equation. Ann. of Math. (2) 137 (1993), no. 2, 295-368.

\bibitem{DGZ03} A. Dembo, A. Guionnet, O. Zeitouni,
Moderate deviations for the spectral measure of certain random
matrices.  Ann. Inst. H. Poincar\'{e} Probab. Statist. 39 (2003),
no. 6, 1013-1042.

\bibitem{DZ98} A. Dembo; O. Zeitouni, Large deviations techniques and applications.
Corrected reprint of the second (1998) edition.
Stochastic Modelling and Applied Probability, 38. Springer-Verlag, Berlin, 2010. xvi+396 pp.

\bibitem{DE13.0}
H. D\"oring, P. Eichelsbacher, Moderate deviations via cumulants. J. Theoret. Probab. 26 (2013), no. 2, 360-385.

\bibitem{DE13}
H. D\"oring, P. Eichelsbacher, Moderate deviations for the eigenvalue counting function
of Wigner matrices. ALEA Lat. Am. J. Probab. Math. Stat. 10 (2013), no. 1, 27-44.

\bibitem{DE13.2}
H. D\"oring,  P. Eichelsbacher,  Edge fluctuations of eigenvalues of Wigner matrices. High dimensional probability VI,
261-275, Progr. Probab., 66, Birkh\"auser/Springer, Basel, 2013.

\bibitem{E12}  L. Erd\"os,
Universality for random matrices and log-gases. Current developments in mathematics 2012, 59-132, Int. Press, Somerville, MA, 2013.

\bibitem{FIK91} A. S. Fokas, A. R. Its, A. V. Kitaev, Discrete Painlev\'e equations and their appearance in quantum gravity.
Comm. Math. Phys. 142 (1991), no. 2, 313-344.

\bibitem{F10}  P. J. Forrester,  Log-gases and random matrices. London Mathematical Society Monographs Series, 34. Princeton University Press, Princeton, NJ, 2010.

\bibitem{G05}  J. Gustavsson,
Gaussian fluctuations of eigenvalues in the GUE.
Ann. Inst. H. Poincar\'{e} Probab. Statist. 41 (2005), no. 2, 151-178.

\bibitem{J98} K. Johansson,
On fluctuations of eigenvalues of random Hermitian matrices.
Duke Math. J. 91 (1998), no. 1, 151-204.

\bibitem{LLX13} S. Z. Li, X. D. Li, Y. X. Xie,  Generalized Dyson Brownian motion, Mckean-Vlasov equation and eigenvalues of random matrices.
arXiv.org/abs/1303.1240v1 (2013).

\bibitem{L16}
D. S. Lubinsky,
Gaussian fluctuations of eigenvalues of random Hermitian matrices associated with fixed and varying weights.
Random Matrices Theory Appl. 5 (2016), no. 3, 1650009, 31 pp.

\bibitem{O09}  S. O'Rourke,
Gaussian fluctuations of eigenvalues in Wigner random matrices.
J. Stat. Phys. 138 (2010), no. 6, 1045-1066.

\bibitem{PS03} L. Pastur, M. Shcherbina,
On the edge universality of the local eigenvalue statistics of matrix models.
Mat. Fiz. Anal. Geom. 10 (2003), no. 3, 335-365.

\bibitem{RS93} L. C. G. Rogers, Z. Shi,
Interacting Brownian particles and the Wigner law.
Probab. Theory Related Fields 95 (1993), no. 4, 555-570.

\bibitem{S09}  M. Shcherbina,
Edge universality for orthogonal ensembles of random matrices.
J. Stat. Phys. 136 (2009), no. 1, 35-50.

\bibitem{S00}  A. Soshnikov, Gaussian fluctuation for the number of particles in Airy, Bessel, sine, and other determinantal random point fields.
J. Statist. Phys. 100 (2000), no. 3-4, 491-522.

\bibitem{S00.2}  A. Soshnikov,  Gaussian limit for determinantal random point fields. Ann. Probab. 30 (2002), no. 1, 171-187.

\bibitem{S06} Z. G. Su,
Gaussian fluctuations in complex sample covariance matrices.
Electron. J. Probab. 11 (2006), no. 48, 1284-1320

\bibitem{TV11}   T. Tao, V. Vu,  Random matrices: universality of local eigenvalue statistics. Acta Math. 206 (2011), no. 1, 127-204.

\bibitem{TW94}  C.  Tracy,  H.  Widom, Level-spacing distributions and the Airy kernel. Comm. Math. Phys. 159 (1994), no. 1, 151-174.

\bibitem{Z14.2} D. Zhang, Random matrices, stochastic nonlinear
Schr\"odinger equations (in Chinese). PhD thesis, Chinese Academy of
Sciences (2014).

\bibitem{Z14}  D. Zhang,  Gaussian fluctuations of eigenvalues in log-gas ensemble:
bulk case I. Acta Math. Sin. (Engl. Ser.) 31 (2015), no. 9, 1487-1500.

\bibitem{Z16}
D. Zhang,
Tridiagonal random matrix: Gaussian fluctuations and deviations.
J. Theoret. Probab. 30 (2017), no. 3, 1076-1103.
\end{thebibliography}
\end{document}